\documentclass[10pt,a4paper]{article}

\usepackage[utf8]{inputenc}
\usepackage{amsmath}
\usepackage{amsfonts}
\usepackage{amssymb}
\usepackage{makeidx}
\usepackage{graphicx}
\graphicspath{ {./} }
\usepackage{listings}
\usepackage{tocloft}
\usepackage{verbatim}
\usepackage{xcolor}
\usepackage{calc}
\usepackage{CJKutf8}
\usepackage[hidelinks]{hyperref}

\usepackage{epigraph}
\newtheorem{le[mm]a}[theorem]{Le[mm]a}

\newcommand{\qed}{\nobreak \ifvmode \relax \else
      \ifdim\lastskip<1.5em \hskip-\lastskip
      \hskip1.5em plus0em minus0.5em \fi \nobreak
      \vrule height0.75em width0.5em depth0.25em\fi}
      
       \usepackage[polutonikogreek,english]{babel}
    
\usepackage[normalem]{ulem}

\usepackage{tikz, tkz-tab, pgfkeys}
\usetikzlibrary{calc,decorations.pathmorphing,shapes}

    \newcommand*{\tg}[1]{\textgreek{#1}}

  \newcounter{exer}[section]

 \title{Notes on angles and solid angles, 
in relation with  Euler's memoir 
\emph{De mensura angulorum solidorum} }
 \author{Stelios Negrepontis and Athanase Papadopoulos}
 
\begin{document}
\maketitle

\vspace*{\fill} \epigraph{\itshape 
For when I say beauty of form, I am trying to express,
not what most people would understand by the words,
such as the beauty of animals or of paintings,
but I mean, says the dialectics,
the straight line and the circle and the plane and solid
figures formed from these by turning-lathes and rulers and patterns of
angles;
perhaps you understand.
For I assert that the beauty of these is not relative, like that of other
things,
but they are always absolutely beautiful by nature
and have peculiar pleasures in no way subject to comparison with the
pleasures of scratching;
and there are colors that possess beauty and pleasures of this character.
Do you understand?
(Plato, \emph{Philebus} 51c-d, 
In:  Plato in Twelve Volumes, Vol. 9, transl. N. Fowler,
 Harvard University Press, 
1925,
with modification by the authors)}

\vspace*{\fill} \epigraph{\itshape 
Je ne parle pas ici, bien entendu, de cette beaut\'e qui frappe les sens, de
la beaut\'e des qualit\'es et des apparences ; non que j'en fasse fi, loin
de l\`a, mais elle n'a rien \`a faire avec la science ; je veux parler de cette
beaut\'e plus intime qui vient de l'ordre harmonieux des parties, et
qu'une intelligence pure peut saisir.
(Henri Poincar\'e, \emph{Science et m\'ethode},  Flammarion, 1908)}

\vfill\eject

\noindent {\bf Abstract} 

We provide some historical context to the study of solid angles carried out by Euler in his memoir \emph{De mensura angulorum solidorum} (On the measure of solid angles). We extend our study to the general notion of angle (not only solid).
While doing so, we explore some works by 
 Ancient Greek mathematicians and others by Arabs mathematicians of the Middle-Ages as well as some later Western authors
 from the Renaissance. In particular, we review the Pythagorean anthyphairetical perspective on angles which establishes the basis of the important relation between the mathematical notion of angle and the philosophical concept of finitization of the Infinite.    In doing so, we shall show that questions addressed by Euler lead us to questions raised about 2500 years ago. At the same time, we highlight the fact that mathematics in those times is also today's mathematics.  The reader can also see in this study the intermingling between mathematics and philosophy. 
 This paper will appear in the book \emph{Spherical geometry in the Eighteenth Century, I: Euler, Lagrange and Lambert}, ed. R. Caddeo and A. Papadopoulos, Springer, 2026.

\bigskip

\noindent {\bf Keywords}  Leonhard Euler, angle, solid angle, Greek mathematics, Euclid, Proclus, Pythagorean mathematics, Plato, anthyphairesis, foundations of mathematics, finitization of the Infinite, history of angles, history of solid angles, Arabic mathematics, 18th century, philosophy of mathematics, Eudoxus of Cnidus.
  
\bigskip

\noindent {\bf AMS classification} 
 01A20, 01A30, 01A50,  	00A30, 03A05.

\bigskip
 
 \tableofcontents

\section{Preamble: On the mathematical notion of angle}

In this chapter, we provide a historical context to the study of solid angles\index{solid angle!Euler} carried\index{Euler, Leonhard!Solid angle} out by Euler in his memoir \emph{De mensura angulorum solidorum} (on the measure of solid angles). At the same time, we survey some works on the more general notion of angle over a period ranging from Ancient Greece until modern times. While doing so, we discuss the idea of angle in Pythagorean thinking where philosophy was intimately related to mathematics, reporting on an  anthyphairetical perspective on angles and the relation the Pythagoreans\index{Pythagoreans} established between\index{Pythagoreans!finitization of the infinite} this notion and their idea of finitizing the Infinite.\index{finitization of the Infinite}

Where should one start, for a historical discussion on the notion of solid angle? And what about the very notion of angle, say planar angle? No doubt, the notion of angle is part of the heritage of all civilizations: Babylonian, Syriac, Egyptian, Indian, Greek, Chinese, etc., and in fact, there is no question that since man was able to rationalize, there was a need, in one way or another, to compare angles, if not to measure them. For example, we can easily imagine a human mentally evaluating the angle the sun makes above the horizon, in order to get an idea of what time of the day it is.\footnote{We note incidentally that the fact of measuring the angle of the sun with the horizon became at some point part of spherical trigonometry. Let us quote the following sentences by Lambert, from the memoir \emph{Anmerkungen und Zus\"atze zur Trigonometrie} (Notes and comments in trigonometry)  whose English translation is contained in the volume \cite{CPST}. From the introduction: ``[\ldots] For example, I explained how one could find the height of the sun for every moment from the length of a day by the simple addition of three logarithms."  And from \S 76 of the same memoir: ``[\ldots] Suppose that given an angle and its two adjacent sides, we want to find the opposite side. This occurs  when calculating tables of solar altitude,\index{calculating tables of solar altitude} and is therefore a frequent problem."} Likewise, at night, by looking at the angle made over the horizon by figures in the sky like the moon, the stars and the constellations, man was able to know (approximately) what time of the night it is.
But at what epoch did the mathematical notion of angle arise? This question leads us to another one: what is a mathematical notion? To what extent is a notion that is described by circles, lines and numbers a mathematical notion?
The answer, if we refer to Bourbaki,\index{Bourbaki, Nicolas} can be found in the opening sentences of his \emph{\'El\'ements de math\'ematique} (first published in 1939): ``Since the Greeks, mathematics has meant proof; some even doubt that outside mathematics there are proofs in the precise and rigorous sense that this word received from the Greeks". A little further, Bourbaki talks about the axiomatic method and formalization. He writes: ``The \emph{axiomatic method} is, strictly speaking, nothing else than this art of writing texts whose formalization is easy to conceive."
Thus, to answer the question of when did a mathematical notion of angle arise, we shall stick to Bourbaki's view on what is mathematics, putting an accent on formalization, axioms, postulates, theorems and proofs, a point of view that we find reasonable, even if it can be contested.

The earliest known texts in which angles are discussed from the point of view of axioms, postulates, theorems and proofs date back to ancient Greece, and this is why, in this overview, we start there. 
Philosophy is also involved in this discussion, and it addresses other kinds of questions.  Indeed, the Pythagorean's mathematics, and, in their wake, Plato's, were intimately linked to their philosophy, and studying  mathematics  without philosophy is pointless when it comes to  Plato\index{Plato} and the Pythagoreans.\index{Pythagoreans} At the same time, treating such a topic in a few pages is a difficult undertaking. We shall nevertheless try to mention a few essential points of this interweaving between philosophy and mathematics, keeping in mind the notion of angle.

Let us start with a mathematical questions about angles. How do we define an angle? And more fundamentally:  is it necessary to define an angle, or is this a primary notion, an undefined object? This was a legitimate question in Euclid's time and it was certainly addressed then. And what about the notion of solid angle? How can it be defined, and how can it be computed?  If we define a solid angle as an area (as Euler does), then what is area? 
Do we need measure theory to define area? 
We read in Heath's translation of the \emph{Elements}\footnote{For the English translation of the \emph{Elements}, we are using throughout this chapter Heath's edition \cite{Heath2}.}  \cite[Definition 11, Book XI]{Heath2}:  ``A solid angle is the inclination constituted by more
than two lines which meet one another and are not in the
same surface, towards all the lines". Is this a clear definition? 

All this is meant to say that the mathematical notion of angle is complex. This is the topic that we shall address in the following pages.


\section{Introduction}

Euler's memoir \emph{De mensura angulorum solidorum}   \cite{Euler-Mensura-T1} is a study of solid angles\index{solid angle!Euler} that culminates in the computation of the measures of those of the five regular convex polyhedra.\index{regular polyhedron!value of angles} At the beginning of his memoir, Euler recalls that the more familiar planar angles are measured by taking a circle of radius 1 centered at the vertex of this angle and measuring the length of the part of the circle that such an angle subtends. Thus, it is natural, by analogy,  to define the measure of a solid angle to be the area of the part subtended by this solid angle on a sphere of radius 1 centered at the vertex of the solid angle. 
Regarding solid angles,
Euler starts by considering the simplest example of solid angle, namely, one obtained by gluing three Euclidean triangles pairwise along their edges to form a vertex of a tetrahedron. Taking a sphere centered at this vertex and intersecting it with the three triangles that form this vertex, we obtain a spherical triangle. Indeed, the intersection of the sphere with a plane passing through its origin is a great circle, and a spherical triangle is precisely a region of the sphere enclosed by three segments of great circles joining pairwise three distinct points. Considering this sphere as a unit sphere, the measure of the solid angle\index{solid angle!Euler} is\index{Euler, Leonhard!Solid angle} then taken to be the area of this spherical triangle. According to a well-known result attributed to Albert Girard (1595-1632), this area is equal to the angle excess of the spherical triangle,\index{Girard, Albert} that is, the excess over two right angles of the sum of its three angles. However, it is natural to ask, independently of this formula, what is meant by area for figures that are not spherical triangles. Indeed, there are solid angles that are cones over figures that are much more general than Euclidean triangles or polygons, and whether such a notion of area was available as a mathematical notion before the modern developments of measure theory is a legitimate question.

 The Greeks had no formal definition of area,\index{area, as a mathematical notion} even though there were some extremely important methods, developed particularly by Archimedes, for computing length,  area and volume. Let us quote  Heath\index{Heath, Thomas} who,\index{Heath, Thomas} in his \emph{History of Greek mathematics} \cite[p. 20]{Heath-History2}, recalls that Archimedes\index{Archimedes} ``performed in fact what is equivalent to \emph{integration} in finding the areas of a parabolic segment, and of a spiral, the surface and volume of a sphere and a segment of a sphere, and the volumes of any segments of the solids of revolution of the second degree." 
 
Returning to the notion of area in Euler's memoir, let us mention that there is also a questioning about the fact of attributing to Albert Girard\index{Girard, Albert} the statement that the area of a spherical triangle is equal to its angle excess.
Indeed, as we recall in a footnote to the English translation of Euler's memoir which is our main subject here, this formula was known, before Girard, to Thomas Harriot\index{Harriot, Thomas}  (1560-1621). Furthermore, Lagrange mentions in his memoir 
\emph{Solution de quelques probl\`emes relatifs aux triangles sph\'eriques avec une analyse compl\`{e}te de ces triangles} (Solution of some problems relative to spherical triangles with a complete analysis of these triangles), included in English translation in the volume \cite{CPST}, that Girard's proof of this theorem is not satisfying, and that the result should rather be attributed to Bonaventura Cavalieri\index{Cavalieri, Bonaventura} (1598-1647) \cite{Cavalieri}. But to simplify, we shall call this result Girard's formula, as Euler does, and as it is often done in the literature on spherical geometry.
The point we want to make here is that even if one can \emph{define} the ``area" of a spherical triangle to be its angle excess, extending this notion to polygons, and then considering this concept as a special case of a general notion of area leads directly to the question of what is area.

After solid angles\index{solid angle!Euler} that are\index{Euler, Leonhard!Solid angle} cones over Euclidean triangles, Euler considers solid angles that are cones over more general Euclidean polygons. But whereas in the first case he deals with triangles with an arbitrary angle at their common vertex, in the more general case, he considers only triangles glued together whose angles at the common vertex are all equal and such that the dihedral angles that any two adjacent triangles make at this vertex are also all equal. In fact, in the application to the computation of the solid angles\index{solid angle!regular polyhedra} of the regular polyhedra\index{regular polyhedron} that Euler carries out in the later sections of his memoir, only this special case is needed. In any case, the measure of the general solid angle\index{solid angle!Euler} that Euler\index{Euler, Leonhard!Solid angle} considers is again taken to be the area of its intersection with a sphere of radius one centered at the given vertex.  
   
While doing so, Euler establishes several formulae for the angle excess of a spherical triangle in terms of its side lengths. These formulae involve trigonometric functions of the sides and they are spherical analogues of the well-known Heron formula\index{Heron formula!spherical} which gives the area of a Euclidean triangle in terms of its side lengths. Presumably, the formulae given by Euler are the oldest known analogues of Heron's formula in spherical geometry. Lagrange, in his memoir \emph{Solution de quelques probl\`emes relatifs aux triangles sph\'eriques avec une analyse compl\`ete de ces triangles}, published 18 years after Euler's one, establishes again such a formula without mentioning Euler's work on this topic.

    As it is often the case in his mathematical memoirs, Euler gives several proofs of these ``spherical Heron formulae" and he treats in detail special cases. Of particular interest is the case where the three sides of the\index{Heron formula!spherical} triangle are infinitely small, in which case he recovers the familiar euclidean Heron formula, verifying the fact that the sphere with its canonical metric (like any other Riemannian manifold) is locally Euclidean.\footnote{Here, we are allowing ourselves to use modern concepts.}
   
    As we already said, Euler's memoir concludes with explicit numerical values that are approximations of the solid angles\index{solid angle!regular polyhedra} of the five regular polyhedra. Euclid's \emph{Elements},\index{Euclid!\emph{Elements}} which culminate in the classification of the five regular polyhedra, do not include such values.  (In fact, Euclid's \emph{Elements} do not include explicit computations.) 
     
The aim of the present chapter is to give a brief (perhaps partly subjective) historical account of the notion of angle, from the Greeks to Euler.

 The content of the rest of this chapter is the following.

 In \S \ref{s:Pythagoreans}, we start by recalling the importance of the notions of angle and solid angle\index{solid angle!Pythagoreans} in the works of the Pythagoreans,\index{Pythagoreans} Plato and Aristotle, an idea that we shall expand in later sections.

In \S \ref{s:Euclid},  we consider angles and solid angles\index{solid angle!Euclid} in Euclid's \emph{Elements}\index{Euclid!\emph{Elements}} and in the works of some later commentators of the \emph{Elements}.

Sections \ref{s:Pythagoreans-anthy} to  \ref{s:finitanization} are concerned with the philosophico-mathematical Pythago\-rean point of view on angle.
We discuss the relation the Pythagoreans made between, on the one hand the notion of acute and obtuse angle,
 and, on the other hand, their principle of the Infinite, the latter being interpreted as the infinite anthyphairesis\index{anthyphairesis!of the diameter to the side of a square} of the diameter to the side of a square.\footnote{We use the term ``diameter" of a square to denote what it called today its diagonal. The Greek word  (\tg{di'ametron}) translates either by diameter or by diagonal. 
Heath writes in \cite[Vol. I, p. 185]{Heath2}: ``Diameter was the
regular word in Euclid and elsewhere for the diameter of a square, and also
of a parallelogram; diagonal (\tg{diag'wnios}) was a later term, defined by Heron
(Def. 68) as the straight line drawn from an angle to an angle".}
This should not be surprising.
   As a matter of fact, the Pythagoreans were so impressed by their discovery of
incommensurability that they considered that the principles of Finite and Infinite, which constituted the main component of their proof of incommensurability, are at the basis of all their mathematics and philosophy. 
The principle of the Finite, which they saw as the cause of the right angle, is based on the finitization of the infinite anthyphairesis of the diameter to the side of a square.\index{finitization of the Infinite}  
This is the way we understand the claim made by Proclus\index{Proclus of Lycia}   (fifth century AD) that the Pythagorean principle of the Finite was for them the cause of the right angle (Proclus 129.17 \cite{Proclus-In-Euclidem}).\footnote{See Proclus, Commentary to the First Book of Euclid's
\emph{Elements}, usually referred to as \emph{In Euclidem}. The standard reference for the English
translation and comments is Glen Morrow's edition \cite{Proclus-In-Euclidem}.}

  The discussion of angles in Pythagorean thought is an illustration of the fact that their mathematics is inseparable from their philosophy,  and also of the fact that they were capable of sweeping
generalizations from very limited evidence, which is also the force of mathematicians.

 In \S \ref{s:Postulate 4}, based on the writings of Proclus, we give our interpretation of Postulate 4 of the \emph{Elements} which says that all right angles are equal.

Since the present essay on angles is a chapter in a book on spherical geometry, it is natural to have a section on angles in non-Euclidean geometry, even if Euler, in his memoir \emph{De mensura angulorum solidorum}, which is the motivation behind this chapter,  does not mention this kind of angles.
	Thus, in \S  \ref{s:angle-parallel}, we review the relation between the notion of angle and the parallel postulate and we take this opportunity to talk about the attempts made by Euler and Lagrange to prove this postulate.

  Section \ref{s:Elements-solid} is a quick overview of the notion of solid angle in the works of Euclid, \index{solid angle!Euclid} Apollonius,\index{Apollonius of Perga} Heron\index{Heron of Alexandria} and Pappus.\index{Pappus of Alexandria}
  
In \S \ref{s:solid-as-magnitude}, we consider the difficulties that arise in considering angles and solid angles as magnitudes.

Section \ref{s:Arabs} is a short review of the notion of angle and solid angle in the works of some Arab mathematicians of the Middle-Ages, who were the natural heirs of the Greek mathematicians of Antiquity.
 
Section \ref{s:Renaissance} contains a quick overview of angles in the works of some Renaissance authors.

In \S \ref{s:angle-spherical}, we recall the notion of angle in spherical geometry.

We conclude the present chapter, in  \S  \ref{s:conclusion}, by pointing out some notions of angle arising in some modern geometrical works.

      For the convenience of the reader, we recall, in a first  appendix to this paper (\S \ref{s:anthy}), the geometric notion of anthyphairesis.
A second appendix, in \S \ref{s:Eudoxus}, concerns Eudoxus of Cnidus, one of the central figures of Greek mathematics and whose name appears several times in this chapter.  Although this takes us a little away from the main subject of our essay, this excursus on Eudoxus will allow us to talk about one of the main figures of Greek science, and to quote Andr\'e Weil,\index{Weil, Andr\'e} Simone Weil\index{Weil, Simone} and Yuri Manin,\index{Manin, Yuri} three mathematicians, philosophers and historians of mathematics who have integrated the meaning of the Greek heritage into their writings. 
   
 \section{The Pythagoreans, Plato and Aristotle}\label{s:Pythagoreans}
 
The Pythagoreans,\index{Pythagoreans} in the sixth century BC, had some quite precise notions of plane and solid angles,\index{solid angle!Pythagoreans} and they used them in their mathematical discoveries, in particular in their work on regular polygons and the classification of regular polyhedra. Less known is the fact that they also established a link between the notion of angle and the concept of Finite and Infinite, via their anthyphairetical proof of the incommensurability\index{Weil, Andr\'e} of the\index{Plato!\emph{Meno}} diameter to the side\index{Weil, Simone} of a square.\footnote{\label{f:Weil} We quote here a letter sent by Andr\'e Weil to his sister Simone, dated March 28, 1940. (The  letter is sent from Rouen, where Andr\'e Weil was in prison \cite[t. VII, Vol. 1, p. 554-555]{Weil-S}.): ``\emph{How} do you think they [the Greeks] discovered the immeasurable? I remember that we don't have any sure data on this, and that I saw somewhere an ingenious suggestion, with an infinite descent where the thing could be \emph{seen} on a square and its diagonal (numerically this gives the following: We have $\frac{\sqrt{2}+1}{1}= \frac{1}{\sqrt{2}-1}$; but if $\frac{a}{b}=\frac{c}{d}$ and $a,b,c$ have a common measure,  
then the same holds for all the lengths obtained by continuing the progression downwards, which is absurd, since these lengths $OB', OC', OB'', OC''$, etc. become as small as we want)."  
[\emph{Comment} crois-tu qu'ils aient d\'ecouvert l'incommensurable~? Je me souviens qu'on n'a pas de donn\'ees certaines l\`a-dessus, et que j'ai vu quelque part une suggestion ing\'enieuse, avec une descente infinie o\`u la chose se \emph{voyait} sur un carr\'e et sa diagonale (num\'eriquement cela donne ceci : on a $\frac{\sqrt{2}+1}{1}= \frac{1}{\sqrt{2}-1}$ ; or, si $\frac{a}{b}=\frac{c}{d}$ et $a,b,c$ ont une commune mesure,  alors il en est de m\^eme de toutes les longueurs obtenues en continuant la progression vers le bas, ce qui est absurde, puisque ces longueurs $OB', OC', OB'', OC''$, etc. deviennent aussi petites qu'on veut).] 
The method of infinite descent suggests anthyphairesis. Among the Ancients, Plato gives an anthyphairetic proof of the incommensurability in the \emph{Meno} 80-86, 97-98, provided one interprets this passage properly.  Once one realizes that there is such a proof in the \emph{Meno},
then it becomes the most ancient preserved proof of the Pythagorean
incommensurability, older by some 35 years of Aristotle's hint of an
arithmetical Archytas type proof that is usually thought to be the most
ancient. One should note that the \emph{Meno}'s proof is not the original one, because it is extremely
streamlined and it also employs the Logos, a concept that was not available to the Pythagoreans. For a discussion of the anthyphairetic proof of incommensurability and its history, see  the recent article by Farmaki and Negrepontis \cite{Negrep-Turaev}.} 
 The discovery of incommensurability was a
brilliant success of the Pythagoreans,\index{Pythagoreans} and their proof is an example of the fact that Mathematics\index{Pythagoreans!finitization of the infinite} is very much the art of
finitizing the Infinite.\index{finitization of the Infinite} About the latter, we can quote a fragment of Philolaos: 

\begin{quote}\small
[Fragment 1] The nature (\tg{`a f'usis})  in the [ordered] world (\tg{'en t\~wi k'osmwi}) is composed  (\tg{\`arm'oqju}) of infinite entities (\tg{'ape'irwn}) and finitizing entities  (\tg{perain'ontwn}), both the world and everything in it.
\end{quote}
This idea is developed in detail in the paper \cite{Negrep-Turaev} by V. Farmaki and the first author of the present chapter.

It may be useful to recall here that definitions and postulates, and in particular those concerning plane angles, did not start with Euclid.  They were established and already used by Greek geometers before Plato.\index{Plato} The latter,\index{Plato!\emph{Republic}} in the \emph{Republic} 510c2-d3, criticizes the geometers ``who base Geometry 
on arbitrary and unknown hypotheses, definitions and postulates", rather than on intelligible, dialectical\index{dialectical principle} principles\footnote{In our interpretation, the word dialectical refers to the Logos criterion of anthyphairesis, as we explain in the sequel.} from which\index{logos criterion} these definitions and postulates result.  
One specific such criticism concerns precisely the (definition of the) three kinds of angles (\tg{gwni\~wn tritt`a eid\~h}, given in the \emph{Republic} 510c4-5), and in particular what became Postulate 4 in Euclid's \emph{Elements},\index{Euclid!\emph{Elements}} the one connected with the right angle.
Plato does not explain how the three kinds of angles and Postulate 4 might follow from Intelligible Beings,\index{Plato!Intelligible Being} but Proclus,\index{Proclus of Lycia} in a long\index{Intelligible Being} passage of his \emph{Commentary to the first Book of Euclid's} \emph{Elements} (\emph{In Euclidem}, 131,9-134, 7) \cite{Proclus-In-Euclidem}, explains that the Pythagoreans,
\index{Pythagoreans} following their epoch-making discovery of the\index{anthyphairesis} anthyphairesis/continued fraction\footnote{For the convenience of the reader, we recall, in the appendix to the present paper, the definition of the anthyphairesis of two magnitudes.} of the 
diameter to the side of a square, and of the diameter itself (\tg{diam'etrou a>ut\~hs}, 510d8), whose sequence of anthyphairetic quotients is [1,2,2,2,\ldots], based the definitions of the three kinds of angles and of Postulate 4 on the closely related to the Pythagorean principle of the Finite and Infinite which
 essentially describe the two fundamental aspects of this anthyphairesis:\index{Pythagoreans!finitization of the infinite} the infinite anthyphairetic division, and its finitization by means of periodicity.\index{finitization of the Infinite}\index{Pythagoreans} 
We shall discuss more thoroughly this question in Sections \ref{s:Pythagoreans-anthy} to 
  \ref{s:finitanization} below.   In the rest of the present section, we shall talk about definitions and postulates before Euclid's times.

Plato\index{Plato} is a heir of the Pythagoreans. One of the early texts that reached us in which solid angles are studied is his dialogue\index{Plato!\emph{Timaeus}} \emph{Timaeus}  \cite{Plato}, written about half a century before Euclid wrote the \emph{Elements}. In this work, Plato exposes his metaphor on the creation of the world which involves  the construction of four regular polyhedra.\index{regular polyhedron} It is generally admitted that Plato learned about the classification of the regular polyhedra from\index{Theaetetus} Theaetetus (ca. 417-369 BC),  a mathematician who was like him a student of Theodorus\index{Theodorus of Cyrene} of Cyrene (ca. 465-398 BC). Theodorus was certainly  influenced by the Pythagoreans.\index{Pythagoreans} In Plato's\index{Plato} mathematical philosophy, the regular polyhedra are neither intelligible nor sensible, but of a third kind, which Plato calls\index{receptacle} the \emph{receptacle}.  The reader interested in Plato's arithmetico-philosophical use of the regular polyhedra and the relation with his notion of the receptacle is referred to the article \cite{Negrepontis2}. 
 Let us now read Plato on regular polyhedra\index{regular polyhedron} in order to see how solid angles\index{solid angle!Plato} enter in his description:

\begin{quote}\small

[\ldots] And when four equilateral triangles are combined so that three plane angles meet in a point, they form {\bf one solid angle, 
which comes next in order to the most obtuse of the plane angles}. 
And when four such angles are produced, {\bf the first solid figure}\footnote{The tetrahedron (pyramid).} is constructed, which divides the whole of the circumscribed sphere into equal and similar\index{Plato!\emph{Timaeus}} parts (\emph{Timaeus} 54e3-55a4).

{\bf And the second solid}\footnote{The octahedron.} is formed from the same triangles, 
but constructed out of eight equilateral triangles, 
which produce\index{solid angle!Plato} {\bf one solid angle out of four planes}; 
and when six such solid angles have been formed, 
the second body in turn is completed (\emph{Timaeus} 55a4-8).

{\bf And the third solid}\footnote{The icosahedron.} is composed of twice sixty of the elemental triangles conjoined, and of twelve solid angles, 
each contained by five plane equilateral triangles, 
and it has, by its production, twenty equilateral triangular bases\index{Plato!\emph{Timaeus}}
(\emph{Timaeus} 55a8-b3).

Now the first of the elemental triangles ceased acting when it had generated these three solids,  
the substance of the {\bf fourth kind}\footnote{The hexahedron (cube).} being generated by the isosceles triangle. Four of these combined, with their right angles drawn together to the center, produced one equilateral quadrangle; and six such quadrangles, 
when joined together, formed eight solid angles,\index{solid angle!Plato} each composed of three plane right angles; and the shape of the body thus constructed was cubic, having six plane equilateral\index{Plato!\emph{Timaeus}} quadrangular bases (\emph{Timaeus} 55b3-c4).
 \end{quote}

 Euler considers in his memoir \emph{De mensura angulorum solidorum} --- which was the motivation for writing the present chapter --- solid angles formed by joining equal triangles. The solid angles of the regular polyhedra which he aims to measure are of this sort. In the passage we recalled above, Plato in the \emph{Timaeus} highlights the fact that the regular polyhedra are constructed using triangles. 
 
  Aristotle\index{Aristotle}, who was definitely attentive to almost all the mathematical questions of his times, was naturally interested in solid angles.\index{solid angle!Aristotle} Heath, in  Vol. I, p. 177, of his edition of the \emph{Elements}, reports on a long debate that took place among philosophers, starting in Antiquity and continuing until the early 20th century Italian commentators, whose object was to know as to what Aristotelian category\index{Aristotelian category!for angles} an angle belongs.
 Among the questions raised there was whether this notions belongs to the category\index{Aristotelian category} of \emph{quantity} (\emph{quantum},  \tg{pos'on} ),  \emph{quality} (\emph{quale}, \tg{poi'on}), \emph{relation} (\tg{pr'os ti}),
 etc.,  whether angles are \emph{magnitudes}, and if yes,\index{magnitude} whether they are \emph{homogeneous} magnitudes,\index{magnitude} whether we can compare angles that are within a certain class, and if yes, what is this class and what are the tools used in such a comparison. There is also the question of whether we may apply to angles the known operations on numbers and magnitudes (addition, multiplication, etc.), the theory of proportions and other mathematical  operations. All these questions were discussed at length by Aristotle in several treatises, including the \emph{Physics}, the \emph{Metaphysics}, \emph{On the Heavens}, the \emph{Categories}, the \emph{Prior} and the \emph{Second Analytics}, and they remained essential in mathematical thinking during the last two and a half millennia. In his comments on the \emph{Elements}, Heath\index{Heath, Thomas} reviews the points of view of several authors on these questions.  Such interrogations were the subject of argumentations between mathematicians, and we shall mention later in this chapter the opinions of several of them on this matter.

After Aristotle, it is natural to talk about Euclid, and this is what we do in more detail in the next section.   
It is likely that the Pythagoreans'\index{Pythagoreans} interest in the regular polyhedra\index{regular polyhedron} and Plato's\index{Plato!regular polyhedra} absorption of this\index{Plato} theory acted as a motivation for Euclid's complete involvement in this topic and his study of solid angles.\index{solid angle!Euclid}

\section{On plane angles in Euclid's \emph{Elements}}\label{s:Euclid}

In Euclid's \emph{Elements},\index{angle!Euclid} angles\index{Euclid!\emph{Elements}} are introduced in Book I. Right angles are first mentioned in the section containing the postulates.   Postulate 4 reads:  \begin{quote}\small
All right angles are equal.
\end{quote}
 Angles also appear at the level of the definitions. Definition 8 of Book I reads:  \begin{quote}\small
 A plane angle is the inclination to one another of two lines in a plane which meet one another and do not lie in a straight line.
 \end{quote} 
 Definition 9 reads:
 \begin{quote}\small
 And when the lines containing the angle are straight, the angle is called rectilineal.
 \end{quote}

 Heath\index{Heath, Thomas} makes the following remark, just after these definitions: ``It looks as though
Euclid really intended to define a rectilineal angle,\index{angle!rectilineal!Euclid} but on second thoughts,
as a concession to the then common recognition of curvilineal angles, altered
`straight lines' into `lines' and separated the definition into two." The remark is interesting and the modern geometers will be interested to know that Euclid also considered plane angles whose sides are not necessarily straight lines.

In any case, these definitions imply that Euclid considers that there exist more than one species of plane angles. Rectilineal\index{angle!rectilineal!Euclid} angles\index{rectilineal angle} form\index{curvilineal angle} a special class of plane angles, the class where the two curves that contain the angle are straight lines. Other classes include the curvilineal\index{angle!curvilineal!Euclid} and the mixed.\index{angle!mixed!Euclid} Among the latter are the so-called semi-circle angle, also\index{horn-like angle} called\index{angle!contact!Euclid} contact, horn-like\index{angle!horn-like!Euclid} or cornicular\index{angle!cornicular!Euclid} angle, the latter name being used by Proclus.\index{Proclus of Lycia} This is the angle made by a circle and a   straight line tangent to it. We shall say more about this angle now.

  Euclid considers non-rectilineal angles in Proposition 16 of Book III and in its proof. The proposition reads:    
  \begin{quote}\small
  The straight line drawn at right angles to the diameter of a circle from its extremity will fall outside the circle, and into the space between the straight line and the circumference another straight line cannot be interposed; further the angle of the semicircle is greater, and the remaining angle less, than  any acute rectilineal angle.
  \end{quote}
    From the proof of the same proposition, we have:

    \begin{quote}\small
    I say further that the angle of the semicircle contained by the straight line $BA$ and the circumference $CHA$ is greater than any acute rectilineal angle, and the remaining angle contained by the circumference $CHA$ and the straight line AE is less than any acute rectilineal angle.
    \end{quote}
    
    Giving an exact meaning to this statement and some others in Euclid's \emph{Elements}  are among the questions that puzzled  ancient mathematicians. Some of these questions found eventually a satisfactory solution with the invention of infinitesimals.\index{infinitesimals}

 It may be good to recall here that Leibniz,\index{Leibniz, Gottfried Wilhelm} when he introduced infinitesimals\index{infinitesimals!Leibniz} (which he also called differentials), postulated them to be greater than zero but smaller than any positive number, and established rules to manipulate them (addition, multiplication, etc.).  Less well known is that in the period between Euclid and Leibniz, some Arab  mathematicians introduced infinitesimals\index{infinitesimals!in Arabic mathematics} and treated them in a similar way, in their work on angles. This is a subject highlighted\index{angle!as magnitude} in the relatively recent book \emph{Angles and magnitude} \cite{Rashed-angles} by Roshdi Rashed on which we shall say more in \S \ref{s:Arabs}.

    Heath has a long historical comment on Proposition III.16, and it is interesting to read it. He starts by saying that ``this proposition is historically interesting because of the controversies to
which the last part of it gave rise from the 13th to the 17th centuries.
History was here repeating itself, for it is certain that, in Ancient Greece, both
before and after Euclid's times, there had been a great deal of the same sort
of contention about the nature of the `angle of a semicircle' and the
`remaining angle' between the circumference of the semicircle and the
tangent at its extremity [\ldots].
Proclus\index{Proclus of Lycia} has no hesitation in speaking of the `angle of a semicircle'  and
the `horn-like angle' as true\index{angle!horn-like!Euclid} \emph{angles}.\index{horn-like angle} Thus he says that `angles are contained
by a straight line and a circumference in two ways; for they are either
contained by a straight line and a convex circumference, like that of the semi-circle, or by a straight line and a concave circumference; like the \tg{keratoeid'hs}' (p. 127-11-14). `There are \emph{mixed} lines, as spirals, and angles, as the angle of the semicircle and the \tg{keratoeid'hs}' (p. 104-16-18).  The difficulty which
the ancients felt arose from the very fact which Euclid embodies in this
proposition. Since an angle can be divided by a line, it would seem to be a\index{angle!as magnitude} 
magnitude; `but if it is a magnitude, and all homogeneous magnitudes which
are finite have a ratio to one another, then all homogeneous angles, or rather
all those on surfaces, will have a ratio to one another, so that\index{angle!cornicular}  the cornicular\footnote{\label{f:Corni} Whereas Euclid considers cornicular angles made by a curve and a line tangent to it. Proclus, in the passage quoted here, considers angles made by a circle and a curve which is not necessarily tangent to it (he talks about an ``angle made by a straight line and a concave circumference").}
will also have a ratio to the rectilineal. But things which have a ratio to one
another can, if multiplied, exceed one\index{Archimedes axiom} another.\footnote{\label{n:Archimedes} The reference is to the so-called Archimedes axiom which says that given two magnitudes $a$ and $b$, there exists a natural number $n$ such that $na>b$. The axiom is also called Eudoxus--Archimedes axiom,\index{Eudoxus--Archimedes axiom} because Archimedes attributes it to Eudoxus. See \S \ref{s:Eudoxus} of this chapter.} Therefore the \emph{cornicular}
angle\footnote{Here, Proclus is talking about a cornicular angle formed by a circle and a line tangent to it.} will also sometime exceed the rectilineal; which is impossible, for it is
proved that the former is less than any rectilineal angle' (Proclus, p. 121,
24-122, 6)."

Furthermore, dealing with angles, compared to lines, needs additional care, since the value of an angle lies within an interval: in Euclid's \emph{Elements}, this value is always assumed to be between $0$ and $\pi$.
How does one then add two arbitrary  angles?

 Thus, in many ways, angles do not satisfy the  Archimedean axiom.\index{Archimedes axiom} This was already stressed by ancient authors. Furthermore, some operations that can be performed on angles cannot be made on angles that do not belong to the same class (recall that Euclid also considers angles that are not rectilineal, for instance the contact angle between a circle and a tangent, an example of a mixed angle).\index{mixed angle}\index{angle!mixed}
Hence the need for a careful discussion of angles. 

It is interesting to note that even though Euclid introduced, in Definitions I.8 and I.9,  a general notion of angle between lines that are not necessarily straight,  he uses the general definition only in the case of the space contained by the straight line drawn at right angles to the diameter
of a circle and the circle itself (Proposition 16 of Book III, which we already mentioned). Furthermore, even though angles bounded by non-straight lines are considered in this very special case, Euclid, in several propositions, is careful in stressing the fact that the angle he refers to is rectilineal\index{rectilineal angle} 
(i.e. bounded by two\index{angle!rectilineal!Euclid} straight lines); see e.g. Proposition I.10 which we recall next, in which Euclid talks about angles made by a ``straight line standing on a straight line"; see also Propositions I.42, I.44, I.45, and there are others.

     The next three definitions of Book I are concerned with the attributes right, acute and obtuse for angles. Definition 10 reads:    
     \begin{quote}\small
     When a straight line standing on a straight line makes the adjacent angles equal to one another, each of the equal angles is right, and the straight line standing on the other is called a perpendicular to that on which it stands.
     \end{quote}
        Definition 11 reads:    
        \begin{quote}\small
        An obtuse angle is an angle greater than a right angle. 
        \end{quote}   Definition 12  reads:     
            \begin{quote} \small
            An acute angle is an angle less than a right angle.
           \end{quote}

Let us turn back now to the Pythagoreans\index{Pythagoreans}  for a more detailed account of their interpretation of the quality for an angle to be acute, obtuse or right.   
    
    \section{The Pythagorean principle of the Infinite, interpreted as the infinite anthyphairesis\index{anthyphairesis!of the diameter to the side of a square} of  the diameter to the side of the square, is the cause of the acute and obtuse angles} \label{s:Pythagoreans-anthy}

    We mentioned in the Introduction that the\index{Pythagoreans} Pythagoreans\index{Pythagorean!principle of the Infinite} established a connection between the principle of the Infinite and the notion of obtuse and acute angles.\index{anthyphairesis} In this section, we shall explain this in some detail, assuming that the reader is less familiar with this kind of questions than with the other purely mathematical questions discussed in the present chapter.

    The key to this connection is the double sequence of the side $p_n$ and diameter $q_n$ of the square, defined by Theon Smyrneus\index{Theon of Smyrna} \cite{Theon}, passages 43,5-15, 44,18-45,8, Iamblichus\index{Iamblichus} \cite{Iamblichus}, passages 21-92,3, 92,23-93,6] and Proclus\index{Proclus of Lycia} \cite{Proclus}, passages 2,24,16-25,13, 2,27,1-29,4, 
     by the recursive rule

\[p_1=q_1=1, \ p_{n+1}=p_n+q_n, \ q_{n+1}=2p_n+q_n \  \mathrm{for} \ n= 1,2,\ldots .\]

 It is easy to prove that this recursive definition is equivalent to the anthyphairetic definition:
 
\[\mathrm{Anth}(q_n, p_n) = [1, 2, 2, 2, \ldots, 2 \ (n-1\ \mathrm{times})] .\]

The line pairs $(q_n, p_n)$ are relatively prime for $n = 1, 2, \ldots $, and their\index{anthyphairesis}  anthyphairesis is a finite initial segment approximating the infinite anthyphairesis of the diameter $d$ to the side $s$ of a square: 

\[\mathrm{Anth}(d,s) = [1, 2, 2, 2, \ldots].\]

For each positive integer $n$, we denote by $\omega_n$ the angle formed by the two equal sides in the isosceles triangle with sides $p_n, p_n, q_n$.
Since the double sequence $(p_n, q_n)$ approximates the dyad $(s,d)$ where $s$ is the side and $d$ the diameter of a square, it is intuitively clear that the angle $\omega_n$ approximates the right angle, without ever being actually equal to it. It is important to realize that this approximation takes the form alternately of an acute and  obtuse angle, exactly because of the fundamental Pell property of the side and diameter numbers, expressed in the next proposition:

\subsection*{Proposition 1} \emph{The sequence of the side and diameter numbers satisfies the Pell property\index{Pell property} 
$q_{2n} = 2p_{2n}^2+ (-1)^n $.}

The proof is by induction, using Proposition II.10 of the \emph{Elements} stating that for any two line segments $a$ and $b$, the equality $(a+2b)^2+a^2=2(a+b)^2+2b^2$ holds.\footnote{The use of the Pell property of the side and diameter numbers provides a strong independent argument that the Pythagoreans\index{Pythagoreans}  had indeed a proof of the Pell property\index{Pell property} of the side and diameter numbers, since now the Pythagoreans actually use this property in order to define the two kinds of angles from the higher principle of the Infinite.\index{Pythagorean!principle of the Infinite} It is safe to assume that the Pythagoreans\index{Pythagoreans}  would not use the Pell property if they had not proved it rigorously. We take this opportunity to note that the use of the Pell property in this Proclus passage shows that, despite the claims of Mueller \cite{Mueller}, Unguru \cite{Unguru} and Fowler \cite{Fowler}, the Pythagoreans not only had knowledge of induction, but in fact employed it\index{Pythagoreans!induction} in order to prove the Pell property of the side and diameter numbers.}

    The Pell property combined with the converse of the Pythagorean theorem (Propositions II.12 and 13 of Euclid's\index{Euclid!\emph{Elements}} \emph{Elements}) has immediate consequences on the kind of angle $\omega_n$, as we shall see now.
    
   \medskip

\subsection*{Proposition 2} \emph{For every $n = 1, 2, \ldots$, the angle $\omega_{2n-1}$ is acute, and the angle $\omega_{2n}$ is obtuse.}

\medskip

\noindent {\it Proof}.  By the Pell property, we have $q_{2n-1}^2 = 2p_{2n-1}^2- 1$, hence $q_{2n-1}^2 < 2p_{2n-1}^2$, therefore by Propositions II.4/5,\footnote{\label{Q45} The numbering II.4/5 indicates the position of the statement and proof of the
Pythagorean theorem in Book II of the
\emph{Elements}, right after II.4 and before II.5, as restored in the paper \cite{Negrepontis3}.} II.12 and II.13, the angle $\omega_{2n-1}$ contained by the equal sides $p_{2n-1}, p_{2n-1}$ of the isosceles triangle with sides $p_{2n-1}, p_{2n-1}, q_{2n-1}$ is acute, and
by the Pell property, $q_{2n}^2 = 2p_{2n+ 1}^2$, hence, $q_{2n}^2 > 2p_{2n}^2$. Therefore, again by Propositions II.4/5, II.12 and II.13, the angle $\omega_{2n}$ contained by the equal sides $p_{2n}, p_{2n}$ in the isosceles triangle with sides $p_{2n}, p_{2n}, q_{2n}$ is obtuse.

Combining the fact that the side and diameter numbers converge to the
dyad (diameter, side) of a square and that this convergence is alternatively
from below and from above, it follows that
the acute angle  $\omega_{2n-1}$  becomes, for sufficiently large $n$, as close to the
right angle as desired, and thus greater than any acute angle,
while the obtuse angle $\omega_{2n}$ becomes, for sufficiently large $n$, as close to
the right angle as desired, and thus smaller than any obtuse angle, thus
finally obtaining the following:
 
\subsection*{Pythagorean definition of the acute and the obtuse angle in terms of the principle of the Infinite}
 
The last sentence of the preceding paragraph justifies the following (Pythagorean) definition of acute and obtuse angle:

 \begin{quote}
\emph{An angle $\omega$\index{Pythagorean!principle of the Infinite} is acute if there is an odd number $2n-1$ such that $\omega < \omega_{2n-1}$, and
an angle $\omega$ is obtuse if there is an even number $2n$ such that $\omega > \omega_{2n}$.}

\end{quote}
This Pythagorean definition of obtuse and acute angle is adopted by
Plato, as evidenced by the\index{Plato!\emph{Republic}} \emph{Republic} passage 510c2-d3, and it becomes a part
of his argument that the geometers are mistaken in relying, for the
foundation of Geometry, on unknown and arbitrary hypotheses (that is, postulates
and definitions), but that they should rather generate Geometry  dialectically,
namely, from his Intelligible Beings,\index{Plato!Intelligible Being} as we already recalled in \S \ref{s:Pythagoreans}.  For this purpose, Plato
repeatedly treats the geometric entity ``diameter to the side of a square",
as an Intelligible Being,\index{Intelligible Being} calling it ``the diameter itself" in the \emph{Republic}
(510d8), and showing,\index{Plato!\emph{Meno}} in the \emph{Meno},\index{Plato!Intelligible Being} 80d-86e and 97a-98b, that its knowledge has the form of ``True
Opinion plus Logos", precisely as the knowledge of an Intelligible Being. See the article \cite{Negrepontis-Meno} for more details about this complex question.

\section{Two more useful arithmetical propositions}\label{s:Propositions}
 
The next two elementary propositions, based on Proposition I.44 of the \emph{Elements}, will be useful in the next two sections.
\medskip

\subsection*{Proposition 3}

 \emph{If $a, b, c, d, e$ are line segments and $A, B, C$ numbers, and if $Aa^2 = Bab + Cb^2$, and $Ac^2 = Bcd + Cd^2$, then $ad = bc$.}

\medskip

\noindent {\it Proof}. By Proposition I.44  of the \emph{Elements}, there is a line segment $e$ such that $ae = bc$.
By assumption, $Aa^2c = Babc + Cb^2c$, 

\noindent whence $Aa^2c = Baae + Cbae$, 

\noindent whence $Aac = Bae + Cbe$, 

\noindent whence $Aacc = Baec + Cbec$, 

\noindent whence $Ac^2a = Bcea + Ce^2a$, 

\noindent  whence $Ac^2 = Bce + Ce^2$. 

Clearly then $d = e$, hence $ad = bc$.

\medskip

\subsection*{Proposition 4} 
\emph{If $a, b, c, d$ are line segments such that $ad = bc$, then $\mathrm{Anth}(a,b) = \mathrm{Anth}(c,d)$.} 
 
 \medskip 
\noindent {\it Proof}. Let $\mathrm{Anth}(a,b) = [k_0, k_1, . . .].$ Thus $a = k_0b + e_1$ with $e_1 < b$. 
Therefore $cb = ad = k_0bd+e_1d$. By Proposition I.44 of the \emph{Elements}, there is a line segment $f_1$ such that $e_1d = bf_1$. Since $e_1 < b$, it is clear that $f_1 < d$. Then $cb = k_0bd+e_1d = k_0bd+bf_1$; therefore $c = k_0d+f_1$ with $f_1 < d$. 

Thus, we have \[\mathrm{Anth}(a,b) = [k_0, k_1,\ldots] = [k_0, \mathrm{Anth}(b,e_1)],\] while \[\mathrm{Anth}(c,d) = [k_0, \mathrm{Anth}(d,f_1)]\] and at this point we have established for the tetrad $b, e_1, d,f_1$ the same condition $bf_1 = de_1$ that we had assumed for the original tetrad $a,b,c,d$ and we are to prove that $\mathrm{Anth}(b,e_1) = \mathrm{Anth}(d,f_1)$. We continue as in the first step, and finish the proof by induction.

\section{The Pythagorean principle of the Finite, interpreted as the Finitization of the infinite anthyphairesis\index{anthyphairesis!of the diameter to the side of a square} of the diameter to the side of a square, is the cause of the right angle}\label{s:finitanization}

 The Pythagorean\index{anthyphairesis}  principle\index{finitization of the Infinite} of the\index{Pythagorean!principle of the Finite} Finite has been\index{Preservation of the Gnomons} interpreted 
as the Preservation of the Gnomons\footnote{The equality $b^2 = 2bc_1 + c_1^2$ can be interpreted as 
the equality of the small square $b^2$ with the Gnomon $2bc_1+c_1^2$ in the greater square $(b+c_1)^2$, hence the (classical) terminology ``Preservation of the Gnomons", an expression well understood by the Pythagoreans, see \cite[Vol. I p. 47]{Heath2}. 
It is probably useful to add a few words about the notion of a gnomon. The Pythagoreans had observed that by adding to a square a figure that has the form of a stalk made up of two perpendicular rectangles of the same thickness, the length of each of them being equal to that of the side of the original square and that are  joined together by a smaller square, one obtains then a new square.
There are other figures that have similar properties. For instance, parallelograms appear  as gnomons in several proofs of propositions in Euclid's \emph{Elements}.  Aristotle\index{Aristotle}, in
 the \emph{Categories} 15a30-31, writes about Gnomons: ``[\ldots] But there are things that increase without altering,
as the square
by surrounding it with a gnomon
is increased but is not thereby altered." In the \emph{Physics} 203a10-16, he writes:  ``Further, the Pythagoreans identify the infinite with the even. For this, they say, when it is cut off and shut in by the odd, provides things with the element of infinity. An indication of this is what happens with numbers. If the gnomons are placed round the one, and without the one, in the one construction the figure that results is always different, in the other it is always the same. But Plato has two infinites, the Great and the Small."
} in the anthyphairesis\index{anthyphairesis}  of the\index{Pythagoreans!finitization of the infinite}  diameter to the side of a square, 
$$b^2 = 2bc_1 + c_1^2,$$

$$c_1^2 = 2c_1c_2 + c_2^2.$$ 

The Preservation of the Gnomons\index{Preservation of the Gnomons} implies, by Proposition 1 of \S \ref{s:Pythagoreans-anthy}, $bc_2 = c_{1}^2$. 

  The condition $bc_2 = c_1^2$, by Proposition 5.2 of the \emph{Elements}, implies that $\mathrm{Anth}(b,c_1)=\mathrm{Anth}(c_1,c_2)$, from which we immediately obtain that $\mathrm{Anth}(a,b)=[1, 2,2,2,\ldots]$.
Thus the Preservation of the Gnomons indeed acts as a Finitizing principle, finitizing our complete knowledge of the infinite anthyphairesis of the diameter to the side of a square.\index{anthyphairesis!of the diameter to the side of a square} 

 On the other hand, from the Preservation of the Gnomons and its consequence   $bc_2 = c_1^2$, since 
$c_1 = a-b$, $c_2 = b-2c_1 = b-2(a-b) = 3b - 2a,$ 
we have, by Proposition II.7,
$$bc_2 =b(3b - 2a) = (a - b)^2 = a^2 + b^2 - 2ab$$
 and
$$c_1^2 =3b^2 - 2ab = a^2 + b^2 - 2ab,$$ 
resulting in  the equality $2b^2 = a^2$. 

Thus, the Preservation of the Gnomons gives the equality $a^2=2b^2$, 
which, by the Pythagorean Theorem I.47 of the \emph{Elements} (and II.4/5 in its restored position) and its converse II.12 and 13, is equivalent to the right angle. In this way,\index{Proclus of Lycia} Proclus' claim that the Pythagorean principle of the Finite is the cause of the right angle is explained.\footnote{ 
Morrow, in \cite[p.104]{Proclus-In-Euclidem}, writes in his Note 93 to Proclus 129.17: 
``A cryptic reference to the distinction between right angles and obtuse or acute angles. See 131.13-134.7."
This remark shows that Morrow has not understood the connection between the Pythagorean treatment of the three kinds of angles and the anthyphairetic proof of incommensurability.}

\section{Our interpretation of Postulate 4 of the \emph{Elements}: All right angles are equal}\label{s:Postulate 4}

Postulate Four of the \emph{Elements} says (see \cite[Vol. I, p. 200]{Heath2}):

\emph{That all right angles are equal to one another}.

It is remarkable that Proclus,\index{Proclus of Lycia} in \emph{In Euclidem}  188,20-189, 12, produces a proof of this postulate, which is perfectly valid in fact, and agreeing in essence with the proof given\index{Hilbert, David} in Hilbert's\index{Hilbert, David!axiomatization of geometry} axiomatisation of Euclidean Geometry, concluding with the comment
 that ``this proof has been given by other commentators and required no great study."

It is then clear that the real content of the Fourth Postulate is in question. 
Proclus in his comments in \emph{In Euclidem} 191, 5-15 suggests that this content is essentially philosophical, strengthening the dialectic status of the right angle:

\begin{quote}\small
This postulate also shows that the rightness of the angles is akin to the equality, 
as the acuteness and obtuseness are akin to the inequality. 
The rightness is in the same column with the equality, 
for both of them belong under the Finite \ldots 
But the acuteness and obtuseness are akin to the inequality \ldots ; for all of them are the offspring (\tg{>'ekgonoi}) of the Infinite. (Proclus, \emph{In Euclidem} 191, 5-11.)
\end{quote}

In view of the anthyphairetic  Pythagorean/Platonic analysis that we pointed out, our interpretation of Postulate 4 of the \emph{Elements}  (all right angles are equal) becomes: 
 \medskip
 
\subsection*{Proposition 5}  \emph{The anthyphairesis of the hypotenuse\index{anthyphairesis}  to the side of every isosceles right  triangle is equal to} [1,2,2,2,\ldots].

 \medskip

\noindent {\it Proof.} Let two isosceles right triangles,
one with right angle $\omega_1$ and with hypotenuse $a$ and sides $b$, and 
the other with right angle $\omega_2$  and with hypotenuse $c$ and sides $d$.
By the Pythagorean theorem I.47 of the \emph{Elements}, $a^2 = 2b^2$,  $c^2 = 2d^2$. 
Then, from the Pythagorean  Principle of the Finite,\index{Pythagorean!principle of the Finite} we obtain, by Proposition 4.1, $ad=bc$, and, by Proposition 4.2, Anth$(a,b)$ = Anth$(c,d)$.  

Thus, by the equality of all right angles (Postulate 4), the Pythagoreans\index{Pythagoreans}  mean that 
in every \emph{isosceles right triangle}, the anthyphairesis\index{anthyphairesis}  of the hypotenuse to the side is equal to 
[1,2,2,2,\ldots]. The postulate is thus seen as a form of Thales' Theorem without ratios of magnitudes.
For more details, we refer the reader to the articles \cite{Negrepontis1} and \cite{Negrepontis3}.

\section{Angles, the parallel postulate, Euler and Lagrange}\label{s:angle-parallel}
To conclude our quick review\index{parallel postulate} of plane angles, let us recall that this notion is at the heart of the work done on the parallel postulate and the attempts that were made to prove it, an activity that occupied some of the greatest mathematicians over a period which lasted more than two thousand years. Let us first recall the statement of this postulate, in Heath's edition \cite[Vol. 1, p. 155]{Heath2}:

\begin{quote} \small That, if a straight line falling on two straight lines
make the interior angles on the same side less than two right
angles, the two straight lines, if produced indefinitely, meet
on that side on which are the angles less than the two right
angles.
\end{quote}

It is well known that there\index{parallel postulate!equivalent forms} are several equivalent forms of this postulate that are formulated in terms of angles. For instance, this postulate can be replaced (keeping all the other postulates unchanged) by one of the following:

\begin{enumerate} 
\item
The angle sum in (one or in) any triangle is equal to two right angles.
\item
The angle sum in a triangle is independent of the chosen triangle.

\item
There exists non equal triangles\footnote{That is (using modern terms), non-isometric triangles.}  that have equal angles.
\item
For any point which is interior to an angle whose measure is strictly between $0$ and $180^{\mathrm{o}}$,   we can draw a line passing through this point and intersecting the two sides of the angle.
\end{enumerate}

Talking about the parallel postulate and since we are especially interested here in the work of Euler, let us recall that 
   Euler himself attempted proofs of this postulate.\index{parallel postulate!Euler attempts of proof} Two\index{Academy!Imperial Saint Petersburg Academy of Sciences} of his attempts were recorded by\index{Fuss, Nikolai Ivanovitch} Nikolai Ivanovitch (or  Nikolaus) Fuss\footnote{Nikola\"\i \  Ivanovitch  Fuss (1755-1826) was, first, the secretary and, then, the assistant of Euler, before becoming his dearest friend and the husband of his grand-daughter Albertine Benedikte Philippine Luise,\index{Euler, Albertine Benedikte Philippine Luise} the oldest daughter of his son Johann Albrecht Euler,\index{Euler, Johann Albrecht} who was also a mathematician. Fuss was himself a very talented mathematician. He helped Euler in reading and writing his correspondence and in writing his memoirs, during the last twelve years of the latter's life, during which Euler was almost completely blind. After Euler's death, Fuss was appointed professor at the Saint Petersburg Imperial Academy of Sciences.  From 1800 and until his death in 1826, he was the permanent secretary of the Academy, and it is at this title that he pronounced a brilliant \emph{Eulogy} of his master, at the latter's death, see \cite{Fuss-Eloge}.  Fuss collected the latter's unpublished works (there were almost 200 such works), he completed those who were not complete, and he prepared them for publication.} in two manuscripts titled \emph{Euleri Doctrina
Parallelismi} and \emph{Euleri Elementa Geometriae ex principio similitudinis
deducta}, which were published in the Soviet Union in  1961 and which are analyzed by Pont in \cite[p. 281-282]{Pont}. 
    
   Euler's first attempt is based on the fact that the set of points in a plane that are a given distance from a line is the union of two lines. The second one
is based on the existence of non-congruent homothetic triangles. We know
that both hypotheses are equivalent to the axiom of parallels. 
   
       Finally, let us also recall that Joseph-Louis Lagrange, who is the second main character in the present book,  \index{parallel postulate!Lagrange attempts of proof}  at some point, believed that he has a proof of the parallel postulate. The following episode is told by Jean-Baptiste Biot\footnote{Jean-Baptiste Biot (1774-1862) was a mathematician, physicist and astronomer. He was at the same time member of the 	
Acad\'emie des sciences and of the Acad\'emie fran\c caise.}
 in \cite[p. 263]{Biot}:
       \begin{quote}\small
       One day, Lagrange took out from his pocket a piece of paper which he read to the Academy\index{Academy!Acad\'emie royale des sciences, Paris} and which contained a proof of Euclid's famous \emph{postulatum} relative to the theory of parallels. This proof was based on an evident paralogism,  which appeared so to everybody, and Lagrange probably realized it while he was reading, because as soon as he finished, he put the piece of paper in his pocket and did not talk about it anymore. A moment of universal silence followed, after which we moved on to other topics.      
        \end{quote}
        
     Lagrange's attempts to prove the parallel postulate are also reviewed by Pont in \cite[p. 231-234]{Pont}; see also Barbarin \cite{Barbarin1928} p. 15 and Bonola \cite[p. 52]{Bonola-1} who quotes de Morgan \cite{Morgan}
           
Now we pass to solid angles.
 
 \section{Solid angles}\label{s:Elements-solid}

   Solid angles  (\tg{stere`a gwn'ia})  appear in Definition 11 and Propositions 20 to 26 of Book XI of the \emph{Elements}.\index{solid angle} Later, they are used in the proof of the fact that there are exactly five regular polyhedra\index{regular polyhedron} (Addendum to Proposition 18 of Book XIII).

  Definition 11 of Book XI reads:

\begin{quote}\small

 A solid angle is the inclination constituted by more
than two lines which meet one another and are not in the
same surface, towards all the lines.
 
 Otherwise: A solid angle\index{solid angle} is that which is contained\index{angle!solid!Euclid} by
more than two plane angles which are not in the same plane
and are constructed to one point. 
\end{quote}

 There is a long historical digression by Heath\index{Heath, Thomas} around this definition, in Vol. II, p. 39-42 of his edition of the \emph{Elements} \cite{Heath2}, covering a period ranging from ancient commentators like Proclus,\index{Proclus of Lycia} until Cardano (1501-1576), Clavius (1538 -1612), Vi\`ete (1540-1603) and others  more recent.
   Heath comments on this ``double definition" (indeed it consists of two statements). He writes that it contains an ambiguity, in particular for the fact that in the first statement Euclid talks about surfaces (\tg{>epif'aneia}) whereas in the second one, he talks only about solid angles bounded by planes. Heath adds that this is a typical case of a definition which is taken in part from older \emph{Elements}.

After Euclid, one may talk about Apollonius of Perga\index{Apollonius of Perga} (3rd-2nd c. BC), who also wrote on solid angles,\index{solid angle} see Heath's\index{Heath, Thomas} comments in \cite{Heath2} and the article by Paul Tannery \cite{Tannery-Fragments}.
Apollonius\index{Apollonius of Perga!solid angle} defined a solid angle\index{solid angle} as the ``bringing together of a surface of a solid to one point as a broken line or surface" (quoted in \cite[vol. III, p. 267]{Heath2}), a definition that Proclus\index{Proclus of Lycia} met with scepticism (see \emph{In Euclidem}, 19-125, 3, quoted in \cite[Chapter 1]{Rashed-angles}). 

Heron of Alexandria\index{Heron of Alexandria} (ca.  10-85 AD) considered solid angles defined by joining to a point arbitrary convex surfaces, see   \cite{Heath2} and \cite{Tannery-Fragments}. Definition 24 of his \emph{Book of Definitions} \cite[vol. IV, p. 28]{Heron-Definitions} says:\footnote{Heath's translation of Euclid's \emph{Elements} Vol. III, p. 267.} ``A solid angle\index{solid angle} is in general the bringing together of a surface which has its concavity in one and the same direction."
  
    Pappus of Alexandria\index{Pappus of Alexandria} (4th c. AD) used solid angles in his work on spherical geometry, namely, in the proof of the spherical triangle inequality,\index{triangle inequality!spherical} given as Proposition 1 of Book VI of his  \emph{Mathematical Collection} \cite[p. 37 of t. II of Ver Eecke's French translation]{Pappus-Collection}. Like Euler, Pappus considered a solid angle\index{solid angle} obtained by gluing three triangles and intersecting it with a sphere centered at the vertex of this angle.  Then, from the triangle inequality
satisfied by the three dihedral angles that form this solid angle  (Proposition XI.20 of the \emph{Elements}), he deduced the triangle inequality for spherical triangles. We note incidentally that Menelaus of Alexandria\index{Menelaus of Alexandria} (ca. 70-140) gave in his \emph{Spherics} a proof of the spherical triangle inequality which is, unlike the one of Pappus, intrinsic to the sphere, that is, that does not make use of the ambient space (this is Proposition 5 of Menelaus' \emph{Spherics},  see \cite[p. 151]{Rashed-Menelaus}).
 
\section{Angles and solid angles as magnitudes}\label{s:solid-as-magnitude}

We already said\index{solid angle!as magnitude} that a discussion of the notion of angle and solid angle in the \emph{Elements} must involve that of magnitude. In fact, a question that arises, since Aristotle,\index{Aristotle} is whether solid angles are magnitudes, and if yes, how, since they do not satisfy the usual supposed properties of a magnitude.
   
Magnitudes\index{magnitude} are not explicitly defined in the \emph{Elements}
 but we are said there that lines (or line segments) are examples of magnitudes: they can be compared and the theory of proportions applies to them. They satisfy a certain number of axioms there. For instance, Axiom 1 of Book I concerns magnitudes; it says:  
\begin{quote}\small
Things which are equal to the same thing are also equal to one another
\end{quote}
 and Axiom 5 of the same book says:

\begin{quote}\small
 The whole is greater than the part.
 \end{quote}

 But can we compare two magnitudes like a straight line and a curve? If yes, how? If this comparison is done by length, how do we compute lengths\index{infinitesimals} without infinitesimal calculus?  Aristotle,\index{Aristotle} in the \emph{Physics} VII and in other writings, already addressed the difficulties encountered in comparing an arc of a circle with a straight line.
 In fact, the theory of proportions developed in the \emph{Elements} applies to magnitudes, but there exist several classes of magnitudes, and only magnitudes of the same class can be compared, added, subtracted, or multiplied. For instance, a surface cannot be compared to a line.  The question is then whether there is a class of magnitudes that comprises angles, and what is exactly this class. 

 Book V of the \emph{Elements}, which contains Eudoxus' theory\index{Eudoxus of Cnidus} of proportion, deals with the class of magnitudes in which every two elements have a certain ratio.\footnote{We have included in the second appendix to the present chapter (\S ef{s:Eudoxus}) a short note on Eudoxus.}  
Definition V.4 says: 
\begin{quote}\small
Magnitudes are said to have a ratio to one another which can, when multiplied, exceed one another.
\end{quote} 

The magnitudes that Euclid considers include lines, areas and solids. Circle segments are examples of non-magnitudes.\index{magnitude} What about angles? Rectilineal angles are measured by the lengths of the segments they make on a circle of radius 1. This constitutes an indication of the fact that they are \emph{not} magnitudes. Let us investigate this more.

One may first recall that there are no explicit computations of measures of length, area, or angle in the \emph{Elements}, except for saying that some angle is right, or acute, or obtuse, or for statements like \emph{The sum of the three angles in a triangle is equal to two right angles} (Proposition 32 of Book I), or \emph{Two circles are to each other like the square of their radii} (Proposition 2 of Book XII).
 Quotients, products, etc. of magnitudes are compared, but never computed. The language is geometrical. For instance, Euclid talks about the ``square on the side", and not the ``square of the side". This point of view is in contrast with that of Archimedes,\index{Archimedes} who had a strong inclination for computations. It is well known that the latter computed approximate values for $\pi$, of areas under a parabola, of spherical discs, etc.

 Difficulties about angles as magnitudes arise for instance when one talks about a cornicular angle, that is, the angle between the circle (or semi-circle) and a tangent line, which is considered in Proposition 16 of Book III of the \emph{Elements} and which we recalled above. Heath,\index{Heath, Thomas} in his edition of the \emph{Elements}, discusses the objections  that some Greek mathematicians made on angle being a magnitude.\index{angle!as magnitude}  In particular, in his long comment on  Proposition III.16, he writes, quoting Proclus and others: 
  \begin{quote}\small
  The difficulty which
the ancients felt arose from the very fact which Euclid embodies in this
proposition. Since an angle can be divided by a line, it would seem to be a magnitude; ``but if it is a magnitude, and all homogeneous magnitudes which
are finite have a ratio to one another, then all homogeneous angles, or rather
all those on surfaces, will have a ratio to one another, so that the \emph{cornicular}
will also have a ratio to the rectilineal. But things which have a ratio to one
another can, if multiplied, exceed one another. Therefore the cornicular
angle will also sometime exceed the rectilineal; which is impossible, for it is
proved that the former is less than any rectilineal angle." (Proclus). The nature of contact between straight lines and circles was
also involved in the question, and that this was the subject of controversy
before Euclid's time is clear from the title of a work attributed to Democritus
(ca. 420-400 BC) [\ldots]
\end{quote}
 
 The precise references to Proclus, Democritus and others are given in Heath's passage. 
Heath\index{Heath, Thomas} then talks about works related to this subject by Cardano, Peletier (Peletarius), Clavius,  Vi\`ete, Galileo Galilei, Wallis and others.

These and other questions on angles, e.g., whether they can be defined as the ``inclination" between two curves, or on how one computes angles that are not bounded by straight lines and how one compare them, lead to  infinitesimal mathematics.\index{infinitesimals} From a modern viewpoint, topological analysis is also involved in this kind of investigation. Indeed, since one defines an angle as the boundary of a region enclosed by two curves, one needs to make precise the notion of ``boundary". Other notions that arise in the context of angles  include continuity, convexity,  infinite division, etc. 
Such questions were raised by Aristotle\index{Aristotle} in various treatises, and they became fundamental objects of research at the Renaissance, finding important developments in the works of Galileo Galilei,\index{Galilei, Galileo} Wallis,\index{Wallis, John} Hobbes\index{Hobbes, Thomas} and other scientists and  culminating in the works of Leibniz\index{Leibniz, Gottfried Wilhelm} and Newton.\index{Newton, Isaac} 

Much less known is that, in the period between the ninth and the thirteenth centuries, these questions were thoroughly studied by Arab mathematicians working in cities like Baghdad, Cairo, Damascus, Aleppo and C\'ordoba, and that their mathematical and philosophical development attained an extremely high level there. This is the subject of the next section.

 \section{Angles and solid angles in the Arabic mathematical literature}\label{s:Arabs} 
   
The Arabic tradition of the \emph{Elements}\index{Euclid!\emph{Elements}!Arabic manuscripts} is still a subject about which  relatively little is known, since a large number of existing manuscripts of Arabic translations or commentaries of (parts of) the \emph{Elements} have not been examined by mathematicians up to the present day. One may quote here Edward Kennedy, from the introduction to his article \emph{Mathematical geography} \cite[p. 185]{Kennedy}: ``The historian of the Islamic exact sciences is frequently confronted with an
\emph{embarras de richesse} --- hundreds of manuscript sources which have never
been studied in modern times." What Kennedy says about general islamic science is true in particular for the literature on Euclid's \emph{Elements} that still has not been seriously examined for the reconstruction of the ``original" version --- if one can still talk about an original version. 
Already at the end of the 18th century, the famous French mathematician and historian of mathematics \index{Montucla, Jean-\'Etienne} Jean-\'Etienne Montucla (1725-1799),  in his \emph{History of Mathematics},  wrote that ``the Arabs provided us with a large number of authors of the same class as Proclus and Theon of Alexandria" \cite[Vol. 1, p. 212]{Montucla}. He mentions, on p. 90 of the same volume, that he had learned that the Oxford Library possessed 400 Arabic manuscripts on astronomy. This is only in one university. It is also worth recalling here that many important Greek mathematical works (including writings by the essential Greek mathematicians like Diophantus,\index{Diophantus of Alexandria} Apollonius,\index{Apollonius of Perga} Archimedes\index{Archimedes} and  Menelaus)\index{Menelaus of Alexandria} have not survived in Greek language, and that it is thanks to the translations, revisions and commentaries made by Arab mathematicians of the Middle-Ages that these texts, or their content in various forms (translations, quotations, paraphrases, commentaries, critiques, corrections, etc.) have come down to us. A large number of still unpublished manuscripts exist. It is worth recalling that Fermat,\index{Fermat@de Fermat, Pierre} Descartes,\index{Descartes, Ren\'e} Euler, Huygens and other prominent mathematicians were devoted readers of works of  Diophantus, Archimedes and others. Christiaan Huygens'\index{Huygens, Christiaan} correspondence with Mersenne and others contains episodes that shows us his efforts to obtain an unpublished Arabic manuscript of Apollonius, which he eventually obtained through the Medici and which he tried to decipher himself, see \cite[t. 2]{Huygens-Correspondance}.
This is also reported on in \cite[p. 42-43]{CHJP}.

Let us return to angles.

On the notion of (plane and solid) angle in Arabic mathematical\index{angle!in Arabic manuscripts} manuscripts, a   706-page book was recently published by Roshdi Rashed, titled  \emph{Angles et grandeur: d'Euclide \`a Kam\=al al-D\=\i n al-F\=aris\=\i }  (Angles and magnitude:  From Euclid to Kam\=al al-D\=\i n al-F\=aris\=\i )  \cite{Rashed-angles}. 
This book involves mathematics, history and philosophy.   
 It contains critical editions together with analyses of a certain number of manuscripts that deal with the notion of angle\index{angle!in Arabic manuscripts} written in Arabic, dating from the early period of Arabic scholarship (ninth century) until the end of that period (13th-14th centuries).  The manuscripts\index{angle!as magnitude} edited in this book include, among others,  translations of parts of Euclid's \emph{Elements}\index{Euclid!\emph{Elements}!Arabic manuscripts} as well as texts originally translated from the Greek on the notion of angle, of which the Greek original is lost.   
 We shall review some of the material contained in this book, especially the part on solid angles.\index{solid angle} The rest of the
  present section may be considered as a quick review of Rashed's book and we do this without going into details and explanations, which would have required an entire volume.

%
%
    
According to Rashed, the first known Arab mathematician who wrote on solid angles is A\d{h}mad ibn `Abd al-Jal\=\i l al-Sijz\=\i\index{Sijz\=\i@al-Sijz\=\i} (ca. 945-1020).  
In a book titled \emph{Introduction to the science of geometry}, al-Sijz\=\i \   gave a classification of the different kinds of solid angles.\index{solid angle!classification} He writes  \cite[p. 328]{Rashed-angles}:
``[\ldots] the first [species of solid angle] under a single surface, the second under a surface and a plane, and the third under planes". 
He gives the following examples of each of these species: ``The angle at the apex of a cone of revolution surrounded by the lateral surface of the cone, or that at one of the apexes of the myrobolan-shaped or quince-shaped solid, for the first;  the angle at one of the vertices of a half-cone, or a half-solid in the shape of a myrobolan or a quince, or a portion thereof, if they are cut by a plane passing through their vertex, for the second; and finally, the angle between three planar angles --- a trihedron --- the sum of which is smaller than four right angles, for the third." The reader will notice that we encounter here solid angles of various shapes,  that are much more general than those alluded to  in Euclid's definition recalled in \S \ref{s:Elements-solid}. Al-Sijz\=\i 's memoir ends with the words: ``These are things concerning angle that leave us puzzled, given that some of its states necessarily imply that it is a magnitude and others that it is not."

%
%

    Among the other texts presented and discussed in Chapter I of Rashed's book, we mention an  anonymous manuscript titled \emph{Treatise on the angle} and referred to as the \emph{Lahore manuscript}, concerned with works of Euclid,\index{Euclid} Apollonius,\index{Apollonius of Perga}   the Byzantine scholars  Agh\=anis\index{Agh\=anis} (fifth century) and Simplicius\index{Simplicius}  (sixth century) and containing a wealth of mathematical proofs and technical remarks on the divisibility of various species of angles.

 Chapter II of the book is concerned more especially with the research conducted\index{magnitude} on the notion of magnitude. It contains in particular a discussion on magnitudes that do not satisfy Definition V.4 of the \emph{Elements} or the so-called Eudoxus--Archimedes\index{Eudoxus--Archimedes axiom} axiom,\footnote{See Footnote \ref{n:Archimedes}.} which is extensively used in Book XII of the \emph{Elements}.\index{Euclid!\emph{Elements}!Arabic manuscripts} We mention in particular, from that chapter, works by  Ab\=u `Ali al-\d{H}asan ibn al-\d{H}asan ibn al-Haytham\index{Ibn al-Haytham} (ca. 965-1040),
whose investigations related to the notion of angle are included notably in two essays he wrote on ``the explanation and correction of Euclid's \emph{Elements}", in which he considers questions like existence, classification, nature and homogeneity of angles. 
Rashed also provides a critical edition of  a text by  the algebraist al-Samaw'al ibn\index{Samaw'al@al-Samaw'al} Ya\b{h}y\=a al-Maghrib\=\i \ (d. 1175)  titled \emph{Epistle on the angle of contact}, in which the latter elaborates on the non-homogeneity and the non-comparability of figures, based on the example of the angle of contact. Together with\index{Farisi@al-F\=aris\=\i}  Kam\=al al-D\=\i n al-F\=aris\=\i \ (d. 1319), al-Samaw'al is one of the successors of   Ibn al-Haytham who continued the latter's research on angles of contact.

  Chapter III of Rashed's book has a more philosophical flavor. It contains an epistle addressed by the famous philosopher and physician Ibn S\=\i n\=a\index{Ibn S\=\i n\=a (Avicenna)} (980--1037), known in the Latin World as Avicenna,\index{Avicenna (Ibn S\=\i n\=a)} addressed to another physician and philosopher, Ab\=u\index{Masihi@al-Mas\=\i \d{h}\=\i} Sahl al-Mas\=\i \d{h}\=\i .  Notions like quantity, quality, relation, magnitude,\index{magnitude} figure, limit of angles, and others are considered from  the philosophical point of view. Infinite divisibility of angles is also discussed. Ibn S\=\i n\=a  considers that angles belong to both categories of quality and quantity.

 Chapter IV is concerned with solid angles.\index{solid angle} It contains notably another text by Ibn al-Haytham,\index{Ibn al-Haytham} who developed the theory of solid angles in the context of  his research on isoperimetry and isoepiphany. This theory was motivated by the problem of  approximating the volume of the sphere by volumes of convex polyhedra, which he considered in the context of his infinitesimal\index{infinitesimals} approach to the study of the sphere. Ibn al-Haytham used the work of Archimedes on the sphere, and also conical projections and spherical geometry. In  Ibn al-Haytham's work, solid angles\index{solid angle} are subject to the usual operations that apply to Archimedean magnitudes\index{Archimedean magnitude} including the theory of proportions.

  Let us quote now a passage of some Arabic versions of the  \emph{Elements}\index{Euclid!\emph{Elements}!Arabic manuscripts} edited in Rashed's book \cite{Rashed-angles} that gives an idea of the questions about solid angles that are dealt with. A version of Euclid's definition\index{solid angle!Euclid} of a solid angle by\index{Abhari@Al-Abhar\=\i}  Al-Abhar\=\i \  (d. in 1265) reads \cite[p. 308]{Rashed-angles}:  

\begin{quote}\small
Solid angle: is any body enclosed by a single surface that terminates at the same point, or more than two plane angles joined at a point, all in the same direction from that point, and such that two of these angles are not in the same plane.
\end{quote}
 
   A mathematician will notice that this definition is clearer than the two definitions in Heath's version of the \emph{Elements} that we recalled in \S \ref{s:Elements-solid}. Here, [Euclid] talks about ``a single surface that terminates at the same point" instead of ``two lines which meet one another and are not in the
same surface, towards all the lines".

  Ibn al-Haytham's\index{Ibn al-Haytham} work on infinitesimal mathematics,\index{infinitesimals} which is a continuation of the works of Archimedes and Apollonius, was already highlighted in another publication\index{infinitesimals!in Arabic mathematics} by Rashed, \emph{Infinitesimal mathematics from the ninth to the 11th centuries}  \cite{Rashed-Infinitesimales}, in 5 volumes. We learn from there that Ibn al-Haytham developed\index{Ibn al-Haytham} a theory of infinitesimals which he used  in the setting of angles, in particular in the comparison between a contact angle and a rectilineal angle. He starts by noting that since these two magnitudes are not Archimedean, a new theory is needed in order to deal with them. Among the arguments he introduced are two sequences, an increasing one and a decreasing one, the second one bounding the first one from above, which he used for the  comparison of \index{infinitesimals!in Arabic mathematics} infinitesimals (see p. 101 of \cite{Rashed-angles}).

  In trying to resolve difficulties that appear in the \emph{Elements}\index{Euclid!\emph{Elements}!Arabic manuscripts} and that are inherent to the notion of angle, and in particular\index{angle!solid!Euclid} of solid angle,\index{solid angle!Euclid}
    Ibn al-Haytham\index{Ibn al-Haytham} wrote two treatises, the {\it Explanation of the postulates of the book of Euclid}, of which a fragment survives and is translated and published in \cite[p. 286 ff.]{Rashed-angles} and the \emph{Book on the resolution of doubts relative to the book of Euclid on the \emph{Elements} and the explanation of its meanings} (see \cite[p. 110 ff.]{Rashed-angles} and \cite{Haytham-Resolution}) in which he addressed a certain number of difficulties.  In these works, he developed a new geometry in which the notion of angle and that of motion\footnote{In modern words, this would be called congruence or isometry.} of figures are primitive elements and where the notion of equality of lines or of areas is  based on motion. We recall here that motion was avoided in Euclid's \emph{Elements} following a ban by Aristotle who\index{Aristotle} considered that this notion pertains to physics and not to mathematics.\footnote{Bourbaki writes, concerning motion \cite[p. 28]{Bourbaki-Histoire}: ``[\ldots] l'usage des
d\'eplacements --- notamment dans les `cas d'\'egalit\'e des triangles'
--- longtemps admis comme allant de soi, devait bient\^ot appara\^\i tre
\`a la critique du XIXe  si\`ecle comme reposant aussi sur des axiomes
non formul\'es." (The use of
motion --- especially in ``cases of equality of triangles" --- long accepted as self-evident, was soon to appear
to 19th-century critics to also be based on unformulated
axioms.)} 
  
  Ibn al-Haytham's work on angles was continued by several Arab mathematicians. We mention in particular a commentary\index{Archimedes!\emph{On the sphere and cylinder}!Arabic manuscripts}  on Archimedes' \emph{Sphere and cylinder} written by Na\d{s}\=\i r al-D\=\i n\index{Tusi@al-\d{T}\=us\=\i , Na\d{s}\=\i r al-D\=\i n} al-\d{T}\=us\=\i \ (1201-1274)  in which the latter addresses the question of the comparability of lines, curves  and curvilineal angles. Na\d{s}\=\i r al-D\=\i n used in particular a notion of ``rolling onto each other," while he was comparing lengths of curves (cf. p. 469 ff. of  \cite{Rashed-angles}). Incidentally,\index{Euclid!\emph{Elements}!Arabic manuscripts} al-\d{T}\=us\=\i \ is also the author of a famous commentary on Euclid's \emph{Elements}. Montucla,\index{Montucla, Jean-\'Etienne}  in his \emph{History of Mathematics} \cite{Montucla}, already wrote in the 18th century that the main editor of Euclid in the East was Na\d{s}\=\i r al-D\=\i n al-\d{T}\=us\=\i . He speaks of his ``learned commentary, printed in Arabic in the magnificent Medici printing house" and adds that ``this work, esteemed among those of his nation, would perhaps also be so among us if a more common language had made it accessible to our understanding."
 
Since Montucla mentions the Medici, let us note that the latter were actively involved in the effort to translate and publish mathematical manuscripts in Arabic. The reader may refer to the recent book by S. Dumas Primbault  \cite{Primbault}, in particular Chapter 3, in which the author recounts how in 1578 the Patriarch of Antioch, visiting Rome, presented to Cardinal Ferdinando de' Medici (later Grand Duke of Tuscany), as a gift, a number of Arabic manuscripts, including a translation of some books of Apollonius' \emph{Conics} that were thought to be lost. The manuscripts found their way into the Medici library before a Latin translation was published a few years later, an undertaking involving several mathematicians, including Galileo Galilei,\index{Galilei, Galileo} Antonio Santini,\index{Santini, Antonio} Vincenzio Viviani,\index{} Giovanni Borelli,\index{Borelli, Giovanni Alfonso} and others. We already mentioned that Christiaan Huygens,\index{Huygens, Christiaan} for mathematical reasons (he was highly interested in Apollonius' works), took a keen interest in the arabic manuscripts of the \emph{Conics}; see his correspondence in \cite[t. 2]{Huygens-Correspondance} and the report on this subject made in \cite[p. 42-43]{CHJP}.

We also quote a passage by 
 Ibn al-Haytham, from his epistle \emph{The explanation of postulates}. He writes \cite[p. 296]{Rashed-angles} : ``Euclid says : `A solid angle is one that is surrounded by planar angles, more than two angles, that are not in the same plane and that meet at a single point.' This assertion is the definition of a solid angle surrounded by flat surfaces.  But there may be solid angles\index{solid angle} that do not fall under this definition; these are angles surrounded by spherical, cylindrical or conical surfaces, except that the solid angles\index{solid angle} he has dealt with in this book are those surrounded by plane surfaces only."  

After Ibn al-Haytham's study of solid angles,\index{solid angle} one had to wait for a few centuries for what regards new developments.

\section{Some notes on solid angles\index{solid angle} from the Renaissance}\label{s:Renaissance}
The history of solid angles starting at the Renaissance\index{Euclid!\emph{Elements}!Renaissance}  involves the works of Johannes Kepler\index{Kepler, Johannes} (1571-1630), Girard Desargues\index{Desargues, Girard} (1591-1661), Ren\'e Descartes (1596-1650), Florimond de\index{Beaune@de Beaune, Florimond} Baune\footnote{Florimond de Beaune (1601--1651) was a jurist and a disciple of Descartes. He is the author of Notes to the latter's \emph{G\'eom\'etrie} published by  van Schooten in  Latin translation. De Beaune was an amateur mathematician and he played a certain role in the diffusion of mathematics in the 17th century.} (1601--1651), Jean-Paul de Gua\index{Gua@de Gua de Malves, Jean-Paul}\footnote{Jean-Paul de Gua de Malves (1710-1786),  known also as l'Abb\'e de Gua, was a 
mathematician, Catholic priest, economist, professor of philosophy at the Coll\`ege de France, translator (from English) and initiator of a project that ended up as the \emph{Encyclop\'edie}.\index{Gua@de Gua de Malves, Jean-Paul!Encyclop\'edie} He was an admirer of Descartes and he published a treatise with the significant title \emph{Usage de l'analyse de Descartes pour d\'ecouvrir, sans le secours du calcul diff\'erentiel, les propri\'et\'es des lignes g\'eom\'etriques de tous les ordres} (The use of Descartes' analysis to discover, without the aid of differential calculus, the properties of geometric curves of all orders). Lagrange, in the memoir translated in \cite{CPST} (Chapter 18), writes: 
``The late de Gua already had the idea of making all spherical trigonometry dependent on a single general property of spherical triangles; but the Memoir that he gave on this subject in the 1783 volume of the Academy of Sciences contains computations that are so complicated that they seem to be more adapted to show the drawbacks of his method than to make it accepted."} (1710--1786) and  many others. Let us say a few words on some of these works.

 Johannes Kepler,\index{Kepler, Johannes} in his \emph{Harmony of the world} \cite{Kepler-Harmony}, considers solid angles\index{solid angle} of\index{solid angle!regular polyhedra} regular\index{regular polyhedron} and semi-regular polyhedra.\index{solid angle!semi-regular polyhedra} In Definition V of Book II, he talks about ``the individual
angles of several plane figures making up\index{solid angle!Kepler}  a solid angle".\index{Kepler, Johannes!solid angle} He studies 3-dimensional tilings, a work in which the notion of solid angle plays an important role. He talks there about ``regular or
semi-regular figures fitted together so as to leave no gap between the sides of the figures". He considers solid figures made of two types of figures fitting together, 
 for instance the rhombic dodecahedron,  a figure having two types of angles: eight obtuse
trilinear solid angles\index{solid angle} and six\index{Euclid!\emph{Elements}!Renaissance} acute quadrilinear solid angles. Kepler writes that ``the result is like the way in which bees
construct their honeycomb: the cells are contiguous and the end of each is surrounded by three opposing ends while its sides are surrounded by the sides of
six more cells. Three more cells could be added at the other end to complete the
figure, except that the entrances to the cells must remain open"  \cite[p. 99-100]{Kepler-Harmony}.

We already mentioned Descartes' \emph{Progymnasmata de solidorum elementis} (Exercises on the elements of solids) \cite{Descartes},\index{solid angle!Descartes} a\index{solid angle!Descartes}  treatise in\index{Descartes, Ren\'e!solid angle} which the French philosopher-mathematician  studies\index{solid angle!of convex polyhedra} solid angles of convex polyhedra, a study containing  a result which is a precursor of Euler's combinatorial formula that gives a relation between the number of vertices, edges and faces of an arbitrary convex polyhedron. 

From Mersenne's \emph{Synopsis mathematica} \cite{Mersenne}, we learn that Desargues\index{Desargues, Girard}  wrote a complete treatise on solid angles. 

Florimond de Beaune's\index{solid angle!Baune@de Baune}   \emph{Doctrine de l'angle solide} (Doctrine of the solid angle) \cite{Beaune} contains\index{Beaune@de Beaune, Florimond!solid angle} an exposition of various geometric construction problems concerning solid angles\index{solid angle} with three faces, given six of their elements (dihedral angles or face angles). 

De Gua's\index{Gua@de Gua de Malves, Jean-Paul!solid angle}  work carries the significant title\index{solid angle!Gua@de Gua}  \emph{Diverses mesures, en partie neuves, des aires sph\'eriques et des angles solides, triangulaires et polygones, dont on est suppos\'e conno\^\i tre des \'el\'emens en nombre suffisant, avec des remarques qu'on croit pouvoir contribuer \`a simplifier les int\'egrations de plusieurs \'equations diff\'erentielles \`a inconnues actuellement s\'epar\'ees} (Various measurements, some of them new, of spherical areas and solid, triangular and polygonal angles, of which we are supposed to know a sufficient number of elements, with remarks that we believe will help simplify the integrations of several differential equations with unknowns that are actually separated)
 \cite{Gua3}.
 
      In all the works mentioned except Euler's, solid angles\index{Euler, Leonhard!Solid angle} are studied,\index{solid angle!Euler} but\index{angle!solid!Euclid} the actual value, or  ``measure", of a solid angle\index{solid angle} is not computed. In contrast, this is what Euler does in his memoir \cite{Euler-Mensura-T1}. 
  It is also worth noting that Euler used in an essential way the notion of solid angle in his proof of the so-called Euler formula\index{Euler, Leonhard!formula for convex polyhedra}; see his memoir  \emph{Elementa doctrinae solidorum} \cite{E230}.
  
In 1925, the biologist and mathematician D'Arcy Thompson\index{Thompson, D'Arcy Wentworth} (1860-1948) published a survey on the semi-regular\index{solid angle!semi-regular polyhedra} polyhedra involving the computation of their angles, see \cite{Thompson}.
 
  To conclude this section, we note that the notion of solid angle\index{solid angle} is widely used in optics, thermal radiation, and other fields of physics.

\section{On the notion of angle in spherical geometry}\label{s:angle-spherical}

	A spherical angle is made by the intersection of two lines (great circles) on the sphere.\footnote{In fact, two great circles make two angles at their intersection point, but in what follows, we shall take the smallest one as the angle made by these two lines.} Euclid, in the \emph{Elements}, even though he works with spheres, does not mention spherical angles. A spherical angle is an intrinsic object of the sphere. Euclid does not study intrinsic properties of the sphere.
	The first geometrical treatise in which the notion of spherical angle is central is the \emph{Spherics} by Menelaus of Alexandria\index{Menelaus of Alexandria} (ca. 70-140 AD). The treatise consists of 91 propositions. Angles appear there generally (but not only) as angles of spherical triangles.\footnote{We recall that on the sphere, considered by all authors of Antiquity as a sphere embedded  in three-dimensional Euclidean space, that a \emph{great circle} is the intersection of this sphere with a plane passing through its centre, that a \emph{segment} on this sphere designates a segment of great circle, and that a \emph{spherical triangle} designates the union of three points on this sphere together with three segments joining them. The segments are called the {sides} of the triangle. For a number of (obvious) reasons, it is generally required that the three points in question and the segments joining them are contained in a half-sphere;  we shall assume without further comments that this holds.}
 A spherical triangle has three sides and three vertices.	
The angle at each vertex of a spherical triangle is defined as the dihedral angle formed by the planes defining the two sides that contain it. The definition of the dihedral angle formed by two planes is given in Euclid's \emph{Elements} (Book XI, Definition 6). 
	
	The first five propositions of Menelaus' \emph{Spherics} are essentially about angles. Proposition 1 concerns the construction of angles, and is used extensively throughout the rest of the treatise. 
	It states that given, on the sphere, a triple consisting of an angle, a line and a point on it, one can construct an angle having this given point as its vertex, a side on the given line and which is equal to the given angle.

	Proposition 5 of the \emph{Spherics} expresses the triangle inequality for spherical triangles.	We have already encountered this proposition in \S \ref{s:Elements-solid}, where we recalled that Pappus used the above definition of angle as a dihedral angle and the fact that 
	the three dihedral angles forming a solid angle satisfy the triangle inequality
 (Proposition XI.20 of the \emph{Elements}) to prove that the triangle inequality is satisfied by the three sides of a spherical triangle. We also mentioned that this is Proposition 1 of Book VI of
Pappus' \emph{Collection}, \cite[p. 37 of t. II of Ver Eecke's French translation]{Pappus-Collection}. Menelaus, in Proposition 5 of his \emph{Spherics}, gives a proof of this spherical triangle\index{spherical geometry!triangle inequality} inequality\index{triangle inequality!spherical} which is, unlike the one of Pappus, intrinsic to the sphere; see \cite[p. 151]{Rashed-Menelaus}.
 
 Let us finally note that in Menelaus' \emph{Spherics}, as in Euclid's \emph{Elements}, angles are compared with each\index{Menelaus of Alexandria} other but their actual measures are never computed. In Menelaus' treatise, there is a definition of acute, obtuse and right angle, and in fact, this is given right at the beginning (first page) of the treatise  (see \cite[p. 504]{Rashed-Menelaus}):
 \begin{quote}\small
 Spherical figures are limited by straight lines, except that their sides are
arcs of great circles, each one of them being smaller than a semicircle. One
which is bounded by three sides is a triangle and, likewise, those which
have four sides. The angles of the triangles are those contained by these
sides.
 
When one of the two planes of the two circles is perpendicular to
another one, then their circumferences intersect with right angles. What is
smaller than that angle is acute and what exceeds it is obtuse.
 \end{quote}

  \section{In guise of a conclusion: Some modern works on the notion of angle}\label{s:conclusion}

   The modern period in mathematics is often considered to be the one starting with Euler.
   Among the ``moderns", we  start by noting that Euler's contemporary, Jean le Rond d'Alembert,\index{Alembert@d'Alembert, Jean le Rond} insisted on clear and simple definitions, accessible to everyone. Naturally, he proposed simple definitions of angles. In his \emph{Essai sur les \'El\'ements de philosophie} (An essay on the Elements of philosophy), he gave a definition of angle that is simple enough but he noted that this definition depends on the notion of line and of the uniformity of the circle, which are not simple notions, see \cite[p. 319]{Alembert-Essai-1}.

  David Hilbert,\index{Hilbert, David} in his \emph{Foundations of geometry} \cite{Hilbert} (first version 1898), a treatise considered by many as an updated  version of the part of Euclid's \emph{Elements} that is concerned with the axioms of geometry, included the notion of angle in the setting of his congruence axioms (Group IV). Felix Klein,\index{Klein, Felix} in his \emph{Elementary mathematics from an advanced standpoint} \cite{Klein} (first ed. 1925) discusses angles at length, making use of the notion of motion, in the third part of this essay, titled \emph{Systematic discussion of geometry and its foundations}.\footnote{This is the part of Klein's treatise that contains a classification of geometries from the point of view of group actions.} In Birkhoff's axiomatization\index{Birkhoff, George David} of geometry \cite{Birkhoff} (1932), which is minimalistic and assumes the real number system, angle belongs to the list of four undefined notions, the other three being point, line and distance.
 
There is a notion of total angle\index{angle!total} at\index{total angle} a point on an arbitrary  surface in the synthetic theory of surfaces developed by A. D. Alexandrov\index{Aleksandrov, A. D.} and the school he founded, see \cite{Alexandrov}. This school is characterized by a return to the techniques of elementary geometry, in the tradition of the Greeks. 
V. A. Zalgaller, in a tribute to Alexandrov on the occasion of the centenary of his birth \cite{Zalgaller}, describing the Leningrad geometry seminar which the latter led in the period 1945-1952, writes: 
``When talking about methodology, [Alexandrov] often pronounced the slogan: `Return to Euclid!', fighting for visual methods."  See also the article by Katateladze  \cite{Kutateladze}, with the significant title 
\emph{Aleksandrov of Ancient Hellas}, in which the author writes that
``Alexandrov enriched the methods of differential geometry by the tools of functional
analysis and measure theory, driving mathematics to its universal status of the
epoch of Euclid." One feature of the modern metric theory of surfaces is that it allows angles of arbitrary value (the value could be greater than four right angles).

Let us also mention the notion of orbifold point,\index{orbifold point} and more generally, of cone points on surfaces, see e.g. Chapter 13 of Thurston's\index{Thurston, William Paul} \emph{Geometry and topology of 3-manifolds} \cite{Thurston-GT3-1}, and his paper \emph{Shapes of polyhedra and triangulations of the sphere} \cite{Thurston-Shapes}. In the latter paper, a model for a  point that has 
 cone angle $\theta\in [0,\infty]$ is the neighborhood of the origin in the metric space obtained
 by taking the universal cover of the Euclidean
plane minus the origin, reinserting the origin and identifying two arbitrary rays issuing from this point by a rotation of angle $\theta$. In the case where $0\leq \theta\leq 2\pi$, such a point with
 cone-angle $\theta$ is also obtained by taking a sector of the Euclidean plane  bounded by two rays emanating from the origin and making an angle $\theta$, and gluing the two boundaries of this sector by a Euclidean rotation of angle $\theta$. A model of the Euclidean cone point of cone\index{cone angle} angle\index{angle!cone} $\theta$ is the neighborhood of the origin in the resulting surface.
 Thurston's theory, like Alexandrov's, is based on a return to the techniques of elementary geometry.
 
Finally, let us note that the notion of solid angle is an essential tool in Sullivan's work on the dynamics\index{Sullivan, Dennis} in hyperbolic space involving measures at infinity, see \cite{Sullivan}.

\section{Appendix: Anthyphairesis} \label{s:anthy}

For the convenience of the reader, we recall the notion of anthyphairesis of two magnitudes, used in several sections of this chapter.

Consider two magnitudes $a$ and $b$, and assume that $a>b$. 
The \emph{anthyphairesis of the pair $(a,b)$} is the sequence of operations
$$a=I_0b+\gamma_1, \ \mathrm{with} \ \gamma_1<b,$$
$$b=I_1\gamma_1+\gamma_2, \ \mathrm{with} \ \gamma_2<\gamma_1,$$
$$\ldots$$
$$\gamma_{n-1}=I_n\gamma_n+\gamma_{n+1}, \ \mathrm{with} \ \gamma_{n+1}<\gamma_n,$$
$$\gamma_{n}=I_{n+1}\gamma_1+\gamma_{n+2}, \ \mathrm{with} \ \gamma_{n+2}<\gamma_{n+1},$$
$$\ldots$$

The resulting sequence of natural numbers $$I_0, I_1, I_2, \ldots$$ is called the \emph{quotients} of the anthyphairesis,
 and $$a>b>\gamma_1 >\gamma_2>\gamma_3>\ldots$$ is the sequence of \emph{remainders} of the anthyphairesis of $a$ to $b$.
 
 The anthyphairesis\index{anthyphairesis} of $a$ to $b$ is \emph{finite}\index{anthyphairesis!finite} if the last remainder divides the immediately previous one, and in this case, this last remainder is the \emph{greatest common divisor} of $a$ to $b$. (This is Propositions VII.1 and 2 of Euclid's \emph{Elements} for natural numbers, and Proposition X.3 for magnitudes.)

 In the case where no remainder divides the immediately previous one, the anthyphairesis\index{anthyphairesis!infinite} is said to be \emph{infinite}. In this case, $a$ and $b$ are \emph{incommensurable} (Proposition X.2 of Euclid's \emph{Elements}).

\section{Appendix: Eudoxus of Cnidus} \label{s:Eudoxus}
 Eudoxus of Cnidus (408-355 BC)  is a mathematician, astronomer, geographer and philosopher to whom we referred in several sections of this chapter. He is an important figure for our subject here because his work has close relations with those the Pythagoreans, of Plato and of Euclid. He was a student of  Archytas\index{Archytas of Tarentum} of Tarentum,  the famous Pythagorean, with whom he studied geometry, number theory and music, among others. Eudoxus visited Plato's\index{Plato's Academy} Academy\index{Academy!Plato} in Athens twice, where he attended Plato's lectures.\footnote{Regarding these two visits, Diogenes Laertius,\index{Diogenes Laertius} in his \emph{Lives of the eminent philosophers} \cite{Diogenes} recounts the following: ``[\ldots] Disembarking at the Piraeus, Eudoxus went up to Athens every day, attended the sophists' lectures, and returned to the port", and then:``He then returned to Athens with a great many disciples; some say he did this to annoy Plato, who had originally snubbed him. Some report that at a drinking party given by Plato, Eudoxus, owing to the number of guests, introduced the practice of arranging couches in a semicircle." Thus, Diogenes clearly implies some friction between Eudoxus and Plato in both visits.
Some modern scholars believe that Plato set Eudoxus as the director of the Academy during his second visit to Sicily, presumably coinciding with Eudoxus' second visit to the Academy. But there is no ancient source supporting this conjecture.}

Eudoxus is considered as the source of a substantial part of Euclid's \emph {Elements}, both for the axiomatic method and for the theory of proportions, the latter being set out in Books V and VI,  and which is at the origin of the modern definition of real numbers. Heath\index{Heath, Thomas} writes in  \cite[p. 326-327]{Heath-History2}: 
``[\ldots] The greatness of the new theory itself needs no further argument when it is remembered that the definition of equal ratios in Eucl. V, Def. 5, corresponds exactly to the modern theory of irrationals due to Dedekind\index{Dedekind, Richard} and that it is word for word the same as Weierstrass's definition\index{Weierstrass, Karl} of equal numbers". 
We already mentioned the Eudoxus--Archimedes axiom,\index{Eudoxus--Archimedes axiom} used in Book V of the \emph{Elements}, which says that, given any two magnitudes $a$ and $b$, there exists a natural number $n$ such that $na>b$ (see Footnote \ref{n:Archimedes}.)

 We also owe to Eudoxus the so-called \emph{exhaustion method},\index{Eudoxus of Cnidus!exhaustion method} also used by Euclid, a theory that represents the beginning of modern integration theory. Bourbaki,\index{Bourbaki, Nicolas} in his \emph{\'El\'ements d'histoire des math\'ematiques}, writes that it would be more appropriate to call the ``Riemann sums'' ``Archimedes sums'', or ``Eudoxus sums'' \cite[p. 248]{Bourbaki-Histoire}. In fact, integration started with
Book XII of the \emph{Elements}, which is due to Eudoxus. It continued and was greatly extended by Archimedes.
Let us quote again Andr\'e Weil,\index{Weil, Andr\'e} from the same letter to his sister Simone\index{Weil, Simone} we quoted in Footnote \ref{f:Weil}. He writes \cite[t. VII, Vol. 1, p. 556]{Weil-S}: ``It seems to me that Eudoxus\index{Eudoxus of Cnidus} should have been the first mathematician in Greek history" [\`A ce qu'il me semble, Eudoxe a d\^u \^etre le premier math\'ematicien dans l'histoire grecque].

On another topic, Eudoxus developed an astronomical theory based on the motions on homocentric  (or concentric) spheres\index{Eudoxus of Cnidus!homocentric spheres} that represent the sun, moon and the five planets known at his epoch (Mercury, Venus, Mars, Jupiter, and Saturn). This system agrees with the known observations of the time.  Eudoxus' theory had a long development throughout the history of Greek mathematical astronomy. It was taken up later and further developed by Aristotle\index{Aristotle} in the \emph{Metaphysics} and the \emph{De Caelo} and it culminated in the work of
Ptolemy\index{Ptolemy of Alexandria} about 500 years later. For a recent commentary on Ptolemy's system\index{Manin, Yuri} of homocentric spheres, see Manin's article \cite{Manin}.
 For two recent 	articles on Eudoxus' works on the theory of proportions (Book V of the \emph{Elements}), on the
method of exhaustion and on his astronomical model of the universe using  homocentric spheres,\index{Eudoxus of Cnidus!astronomical system} see the recent articles \cite{Negrepontis1} by Negrepontis and \cite{Sunada} by Sunada.  Let us also quote Simone Weil, from her work \emph{Intuitions pr\'e-chr\'etiennes/Esquisse d'une histoire de la science grecque} \cite{Weil-Intuitions}.
\begin{quote}\small
The same Eudoxus developed an astronomical system to answer Plato's question: ``To find the set of uniform circular motions that allow us, concerning the stars, to keep up appearances."
It relies on the brilliant idea of composition of motions that is at the basis of our mechanics.  Just as we construct the parabolic movement of projectiles with two rectilineal motions, one uniform, the other accelerated, so Eudoxus accounted for all the positions of the stars observed in his time by assuming several uniform circular motions accomplished simultaneously around different axes by a single star.
This conception is as daring, as pure as that which defines real numbers and that of the summation of an infinite sum. If Plato wanted only uniform circular motions, it is because this motion alone is divine, and because the stars are, he says, images of the divinity sculpted by the divinity itself. Plato almost certainly has this composition of motions in mind when he talks about the Other, rebellious to unity, harmonized by constraint with the Same. The sun, in its unique motion, is carried along at the same time by the circle of the equator and by that of the ecliptic, which correspond to the Same and the Other; and this makes only one motion.  
In the following period of Greek science, Ptolemy reproduced, in a much cruder form, Eudoxus' system; Apollonius continued  Menechmus\footnote{Menaechmus (c. 380 - c. 320 BC) was a geometer and philosopher, disciple of Plato and of Eudoxus. To him is attributed the discovery of conic sections and his solution to the  problem of doubling the cube using the parabola and hyperbola.} discoveries of conics, and Archimedes those of Eudoxus on integration.
\end{quote}

  The ancient sources for Eudoxus include the \emph{Lives} by Diogenes La\"ertius \cite{Diogenes}, a work we already quoted, and passages from Aristotle's \emph{Metaphysics},  \emph{Nicomachean Ethics} and other works, and from Archimedes, Callippus, Proclus, Simplicius,  and there are others.

\end{document}